\documentclass[notitlepage,leqno,10pt]{article}

\usepackage{graphicx,multirow}
\usepackage{amssymb,amsmath,amsthm}
\usepackage{lscape,lipsum,microtype}
\usepackage[margin=2.54cm]{geometry}  
\usepackage{hyperref}

\usepackage{ifluatex}
\ifluatex
 \usepackage{unicode-math}
\fi

\usepackage{latexsym}
\usepackage{amsmath}
\usepackage{amsfonts}
\usepackage{amssymb}
\usepackage{amscd}

\renewcommand{\d}{\delta }
\newcommand{\D }{\Delta }

\renewcommand{\l }{\lambda }
\renewcommand{\L }{\Lambda }
\newcommand{\m }{\mu }
\newcommand{\n }{\nabla }

\newcommand{\Sig }{\Sigma}

\newcommand{\ov}{\overline}
\newcommand{\intbar}{\mathop{\int\makebox(-13.5,0){\rule[4pt]{.7em}{0.3pt}}%
\kern-6pt}\nolimits}

\newcommand{\be}{\begin{equation}}
\newcommand{\ee}{\end{equation}}
\newcommand{\bes}{\begin{equation*}}
\newcommand{\ees}{\end{equation*}}
\newcommand{\ba}{\begin{eqnarray}}
\newcommand{\ea}{\end{eqnarray}}
\newcommand{\bas}{\begin{eqnarray*}}
\newcommand{\eas}{\end{eqnarray*}}
\newenvironment{pf}{\noindent{\sc Proof}.\enspace}{\rule{2mm}{2mm}\medskip}
\newenvironment{pfn}{\noindent{\sc Proof}}{\rule{2mm}{2mm}\medskip}

\newcommand{\R}{\mathbb{R}}

\newcommand{\Z}{\mathbb{Z}}

\newcommand{\N}{\mathbb{N}}

\author{ Cheikh Birahim NDIAYE}

\date{}

\title{\bf Leray-Schauder degree for the resonant \:$Q$-curvature problem in even dimensions}
\begin{document}

\newtheorem{lem}{Lemma}[section]
\newtheorem{pro}[lem]{Proposition}
\newtheorem{thm}[lem]{Theorem}
\newtheorem{rem}[lem]{Remark}
\newtheorem{cor}[lem]{Corollary}
\newtheorem{df}[lem]{Definition}

\maketitle

\begin{center}

{\small

\noindent  Department of Mathematics Howard University \\  Annex $3$, Graduate School of Arts and Sciences, $217$\\ DC 20059 Washington, USA

}
{\small

\noindent

}
\end{center}

\footnotetext[1]{E-mail addresses: cheikh.ndiaye@howard.edu\\
\thanks{\\ The author was partially supported by NSF grant DMS--2000164.}}

\

\

\begin{center}
{\bf Abstract}
\end{center}
In this paper, using the theory of critical points at infinity of Bahri\cite{bah}, we derive an exact bubbling rate formula for the resonant prescribed \;$Q$-curvature equation on closed even-dimensional Riemannian manifolds. Using this, we derive new existence results for the resonant prescribed \;$Q$-curvature problem under a positive mass type assumption. Moreover, we derive a compactness theorem for conformal metrics with prescribed \;$Q$-curvature under a non-degeneracy assumption. Furthermore, combining the bubbling rate formula with the construction of some blowing-up  solutions, we compute the Leray-Schauder degree of the resonant prescribed \;$Q$-curvature equation under a non-degeneracy and Morse type assumption. 
\begin{center}


\bigskip

\centerline{\bf AMS subject classification: 53C21, 35C60, 58J60.}

\end{center}

\section{Introduction and statement of the results} 

One of the most recurrent question in conformal geometry is the problem of finding conformal metrics for which a certain curvature quantity is equal to a prescribed function, e.g. constant. As a model of such problems, we have the problem of existence of conformal metrics with prescribed Gauss curvature on closed Riemannian surfaces, namely the Nirenberg problem and the more general Kazdan-Warner problem. 
\vspace{6pt}

\noindent
There exists also analogues of the Gauss curvature in high even dimensions which enjoy similar properties which are relevant to conformal geometry. To better introduce those curvatures, we recall some facts about the theory of closed Riemannian surfaces. It is a well known fact that the Laplace-Beltrami operator on  closed Riemannian surfaces\;($\Sigma, g$) is conformally covariant of bidegree $(0, 2)$, and  governs the transformation laws of the Gauss curvature under conformal changes of the background metric $g$. In fact, under  the conformal change of metric \;$g_u=e^{2u}g$, we have
\ba\label{eq:g}
\Delta_{g_u}=e^{-2u}\Delta_{g};\;\;\;\;\;\;\;\;\;-\Delta_{g}u+K_{g}=K_{g_u}e^{2u},
\ea
where\;$\Delta_{g}$\;and\;$K_{g}$\;(resp.\;$\Delta_{g_u}$\;and\;$K_{g_u}$)\;are the Laplace-Beltrami operator and the Gauss curvature of ($\Sigma, g$) (resp. of ($\Sigma,  g_u$)). Moreover, we have the Gauss-Bonnet formula which relates \;$\int_{\Sigma}K_{g}dV_{g}$\;and the topology of\;$\Sigma$
$$
\int_{\Sigma}K_{g}dV_{g}=2\pi\chi(\Sigma),
$$
where \;$\chi(\Sigma)$\; is the Euler characteristic of \;$\Sigma$\; and \;$dV_g$ is the volume form on \;$\Sig$\; with respect to \;
$g$. From this we have that \;$\int_{\Sigma}K_{g}dV_{g}$\; is a topological invariant, hence a conformal invariant too. We point out that the conformal invariance of \;$\int_{\Sigma}K_{g}dV_{g}$ can also be seen just by integrating equation \eqref{eq:g} which is of divergence structure.
\vspace{6pt}

\noindent
In 1983, Paneitz\cite{p1} has discovered a conformally covariant differential operator \;$P_g$\; on four-dimensional closed Riemannian manifolds \;$(M, g)$\; (known now as Paneitz operator). To the Paneitz operator, Branson\cite{bran1} has associated a natural curvature invariant \;$Q_g$\; called \;$Q$-curvature. They are defined in terms of the Ricci tensor \;$Ric_{g}$\; and the scalar curvature \;$R_{g}$\; of the Riemannian manifold \;$(M,g)$ as follows
\ba\label{eq:P}
P_{g}=\D_{g}^{2}+div_{g}\left((\frac{2}{3}R_{g}g-2Ric_{g})\n_g\right),\;\;\;Q_{g}=-\frac{1}{12}(\D_{g}R_{g}-R_{g}^{2}+3|Ric_{g}|^{2}),
\ea
\noindent
where \;$div_g$\; is the negative divergence and \;$\n_g$\; is the gradient with respect to \;$g$. As the Laplace-Beltrami operator is conformally covariant of bidegree \;$(0, 2)$, and governs the transformation laws of the Gauss curvature under conformal changes, we have also that the Paneitz operator is conformally covariant of bidegree \;$(0, 4)$, and governs the transformation laws of the \;$Q$-curvature under conformal changes of the background metric. Indeed under a conformal change of metric \;$g_u=e^{2u}g$, we have
\begin{equation}\label{eq:tranlaw}
P_{g_u}=e^{-4u}P_{g},\;\;\;\;\;\;\;\;\;P_{g}u+2Q_{g}=2Q_{g_u}e^{4u}.
\end{equation}
Apart from this analogy, we have also an extension of the Gauss-Bonnet identity which is the  Chern-Gauss-Bonnet formula
$$
\int_{M}(Q_{g}+\frac{|W_{g}|^{2}}{8})dV_{g}=4\pi^{2}\chi(M),
$$
where \;$W_{g}$\; denotes the Weyl tensor of \;$(M,g)$\; and \;$\chi(M)$\; is the Euler charcteristic of \;$M$.\;Hence from the pointwise conformal invariance of \;$|W_{g}|^{2}dV_{g}$,\;it follows that \;$\kappa_{P}:=\int_MQ_gdV_g$\; is also conformally invariant. Moreover, the conformal invariance of \;$\int_{M}Q_{g}dV_g$\; can also be seen by just integrating the second equation of \eqref{eq:tranlaw} which is also of divergence structure, like in the case of Gauss curvature for closed Riemannian surfaces.
\vspace{6pt}

\noindent
On the other hand, there are high order analogues of the Laplace-Beltrami operator and  the Paneitz operator for high even-dimensional closed Riemannian manifolds and also to the associated curvature invariants. More precisely, given a closed \;$n$-dimensional Riemannian manifold \;$(M, \;g)$\; with \;$n\geq 2$\; and even,\;in \;\cite{gjms},\;it was introduced a family of conformally covariant differential operators \;$P^{n}_{g}$\; whose leading term is \;$(-\D_{g})^{\frac{n}{2}}$. Moreover, in\;\cite{bran1},\;some curvature invariants \;$Q_{g}^{n}$\; was defined, naturally associated to\;$P^n_{g}$. In low dimensions, we have the following relations
$$
P^{2}_{g}=-\D_{g},\;\;\;\;\;\;Q^{2}_{g}=K_{g},\;\;\;\;P^{4}_{g}=P_{g},\;\;\;\text{and}\;\;\;Q^{4}_{g}=2Q_{g}.
$$
It turns out that\;$P^n_g$\;is self-adjoint and annihilates constants. Furthermore, as for the Laplace-Beltrami operator on closed Riemannian surfaces and the Paneitz operator on closed four-dimensional Riemannian manifolds, for every closed \;$n$-dimensional Riemannian manifold \;$(M, g)$\; with \;$n\geq 2$ and even, we have that after a conformal change of metric \;$g_u=e^{2u}g$
\begin{equation}\label{eq:confev}
P^{n}_{g_{u}}=e^{-nu}P^{n}_{g},\;\;\;\;\;\;\;\;P^{n}_{g}u+Q^{n}_{g}=Q^{n}_{g_{u}}e^{nu}.
\end{equation}
Thus, as in the \;$2$-dimensional and \;$4$-dimensional cases, by integrating the second equation of \eqref{eq:confev} and using the fact that \;$P^n_g$\; is self-adjoint and annihilates constants, it follows that \;$\kappa_g^{n}:=\int_{M}Q^n_gdV_g$ is also conformally invariant.  Furthermore, there exists also an analogue of the Chern-Gauss-Bonnet formula, which is a consequence of a formula of Alexakis\cite{alex1} (see also \cite{alex2}), and reads as follows
\begin{equation}\label{eq:alexakis}
\int_M(Q^n_g+|\tilde W_g|)dV_g=\frac{(n-1)!}{2}\omega_n\chi(M),
\end{equation}
 where \;$\omega_n$\; is the volume of \;$\mathbb{S}^n$ (the \;$n$-dimensional unit sphere of \;$\R^{n+1}$) with respect to the standard metric of \;$\mathbb{S}^n$,  $\tilde W_g$\; is a local conformal invariant involving the Weyl tensor \;$W_g$\; and its covariant derivatives, and \;$\chi(M)$ is the Euler characteristic of \;$M$.
 \vspace{6pt}
 
\noindent
As for the Kazdan-Warner problem for closed Riemannian surfaces, given a closed \;$n$-dimensional Riemannian manifold \;($M,g$) ($n$\; even and \;$n\geq 4$) and a smooth positive function \;$K: M\longrightarrow R$, one can ask under which conditions on \;$K$, does there exists a Riemannian metric \;$g_u=e^{2u}g$\; conformally related to \;$g$\; for which the corresponding \;$Q$-curvature \;$Q^{n}_{g_u}$\; is \;$K$.  Thanks to \eqref{eq:confev}, the geometric problem is equivalent to finding a smooth solution of the equation
 \begin{equation}\label{eq:qequationtrue}
P_{g}^nu+Q_{g}^n=K e^{nu}\;\;\;\;\;\;in\;\;M.
\end{equation}
Since \;$\R\subset \ker P^n_g$, then it is easy to see that \eqref{eq:qequationtrue} has a solution is equivalent to 
\begin{equation}\label{eq:qequation}
P_{g}^nu+Q_{g}^n=\kappa_g^nK e^{nu}\;\;\;\;\;\;in\;\;M.
\end{equation}
has a solution. Equation \eqref{eq:qequation} has a variational structure. Indeed, using elliptic regularity theory (see \cite{uv}), we have smooth solutions of \eqref{eq:qequation} can be found by looking at critical points of the following geometric functional:
$$
J(u):= \langle P^n_g u,u\rangle \, + \, 2 \int_M Q_g^n u dV_g \,  - \,\frac{2}{n}\kappa_{g}^{n} \log\left(\int_{M}Ke^{nu} dV_g\right),\;\;\;\;u\in W^{\frac{n}{2}, 2}(M),
$$
where \;$W^{\frac{n}{2}, 2}(M)$\; is the space of functions on \;$M$\; which are of class \;$W^{\frac{n}{2}, 2}$\; in each coordinate system.
\vspace{6pt}

\noindent
The asymptotic behaviour of sequences of solutions \begin{equation}\label{eq:blowupeq}
P_g^nu_l+t_lQ_g^n=t_l \kappa_g^n Ke^{nu_l}\;\;\;\text{in}\;\;\;M,
 \end{equation}
with \;$\ker P_g^n\simeq \R$, $K: M\longrightarrow \R$\; smooth and positive, and \;$t_l\longrightarrow 1$\; as $l\longrightarrow +\infty$\; plays an important role in the Variational Analysis of \;$J$\; in the resonant case, ie when \;$\kappa_g^n\in (n-1)!\omega_n\N^*$. In this paper, we are interested in the exact bubbling rate formula for \eqref{eq:blowupeq} in the resonant case and its applications to existence, compactness and topological degree-computation.
\vspace{6pt}

\noindent
In order to state our results clearly, we first fix some notation and make some definitions. For \;$m\in \N^*$\; such that \;$$\kappa_{g}^n=(n-1)!m\omega_n,$$ we define \;$\mathcal{F}_K: M^m\setminus F_m(M)\longrightarrow \R$\; as follows
\begin{equation}\label{eq:limitfs}\\
\mathcal{F}_K(a_1,\cdots, a_m):=\sum_{i=1}^m\left(H(a_i, a_i)+\sum_{j\neq i}G(a_i, a_j)+\frac{2}{n}\log(K(a_i))\right),
\end{equation}
where \;$F(M^m)$ denotes the fat Diagonal of $M^m$, namely $F(M^m):=\{A:=(a_1, \cdots, a_m)\in M^m:\;\;\text{there exists} \;\;i\neq j\;\,\text{with}\;\, a_i=a_j\}$, $G$ is the Green's function of $P_g^n(\cdot)+\frac{1}{m}Q_g^n$\; with mass \;$(n-1)!\omega_n$\; satisfying the  normalization $\int_M Q_g^n(x) G(\cdot, x)dV_g(x)=0$, and $H$ is its regular part, see Section \ref{eq:notpre} for more information. Furthermore, we define
\begin{equation}\label{eq:critfk}
Crit(\mathcal{F}_K):=\{A\in M^m\setminus F_m(M): \;\;A\;\;\;\text{critical point of} \;\;\mathcal{F}_K\}.
\end{equation}
Moreover,  for $A=(a_1,\cdots, a_m)\in M^m\setminus F_m(M)$, we set
\begin{equation}\label{eq:partiallimit}
\mathcal{F}^A_i(x):=e^{n(H(a_i, x)+\sum_{j\neq i}G(a_j, x))+\frac{1}{n}\log(K(x))},
\end{equation}
and define
\begin{equation}\label{eq:defindexa}
\mathcal{L}_K(A):=-\sum_{i=1}^m  (\mathcal{F}^{A}_i)^{\frac{6-n}{2n}}(a_i)L_g((\mathcal{F}^{A}_i)^{\frac{n-2}{2n}})(a_i),
\end{equation}
where $$L_g:=-\D_g+\frac{(n-2)}{4(n-1)}R_g$$ is the conformal Laplacian associated to $g$.  We set also
\begin{equation}\label{eq:critsett}
\mathcal{F}_{\infty}:=\{A\in Crit(\mathcal{F}_K):\:\;\mathcal{L}_K(A)<0\},
\end{equation}
and
\begin{equation}\label{eq:minf}
i_{\infty}(A):=(n+1)m-1-Morse(A, \mathcal{F}_K),
\end{equation}
where \;$Morse(\mathcal{F}_K, A)$\; denotes the Morse index of \;$\mathcal{F}_K$ at  $A$.
\vspace{4pt}
 Finally, we say
\begin{equation}\label{eq:nondeg}
(ND)_0 \;\;\;\text{holds if for every}\;\;A\in Crit(\mathcal{F}_K) , \;\;\mathcal{L}_K(A)\neq 0,
\end{equation}
\begin{equation}\label{eq:nondeg}
(ND)_- \;\;\;\text{holds if for every}\;\;A\in Crit(\mathcal{F}_K) , \;\;\mathcal{L}_K(A)< 0.
\end{equation}
\begin{equation}\label{eq:nondeg}
(ND)_+ \;\;\;\text{holds if for every}\;\;A\in Crit(\mathcal{F}_K),  \;\;\mathcal{L}_K(A)>0,
\end{equation}
and
\begin{equation}\label{eq:nondeg}
(ND) \;\;\;\text{holds if}\; \;(ND)_0 \;\;\text{holds and}\;\;\mathcal{F}_K\;\; \text{is a Morse function}.
\end{equation}
\vspace{6pt}

\noindent
Now, we are ready to state our results and we start with the exact bubbling rate formula for sequences of blowing up solutions to \eqref{eq:blowupeq} under the assumption \;$\ker P_g^n\simeq \R$. 
\begin{thm}\label{eq:esttl}
Let \;$(M, g)$\; be a closed \;$n$-dimensional Riemannian manifold with \;$n\geq 4$\; even such that \;$\ker P_g^n\simeq \R$\; and  \;$\kappa_g^n=(n-1)!m\omega_n$\; with \;$m\in \N^*$. Assuming that \;$K$\; is a smooth positive function on \;$M$\; and \;$u_l$\; is a sequence of blowing up solutions to \eqref{eq:blowupeq} with \;$t_l\rightarrow 1$\; as \;$l\rightarrow +\infty$, then up to a subsequence, we have that for \;$l$\; large enough, there holds
$$
t_l-1=\frac{c_{n, m}^K(A)e^{-2\max_{M}u_l}}{(\mathcal{F}^{A}_i(a_i))^{\frac{n-2}{n}}}\left[\mathcal{L}_K(A)+o_l\left(1\right)\right],
$$
with \;$A=(a_1, \cdots, a_m)\in Crit(\mathcal{F}_K)$\; and \;$c_{n, m}^K(A)$ is a positive constant depending only on $K$, \;$A$, \;$n$\; and \;$m$\;.
\end{thm}

\vspace{6pt}
 
\noindent
Theorem \ref{eq:esttl} and our work\cite{nd1} in the {\em nonresonant} case (i.e \;$ \kappa_g^n\notin(n-1)!\omega\N^*$) imply the following existence result for conformal metrics with prescribed \;$Q$-curvature.
\begin{cor}\label{eq:corexistence}
 Let \;$(M, g)$\; be a closed $n$-dimensional Riemannian manifold with \;$n\geq 4$\;even such that \;$\ker P_g^n\simeq \R$\; and \;$\kappa_g^n=(n-1)!m\omega_n$\; with \;$m\in \N^*$. Assuming that \;$K$\; is a smooth positive function on \;$M$\; such that $(ND)_{-}$ or $(ND)_+$ holds, then \;$K$\; is the \;$Q$-curvature of a Riemannian metric conformally related to $g$.
\end{cor}

\vspace{6pt}

\noindent
In the critical case, i.e \;$\kappa_g^n=(n-1)!\omega_n$, Theorem \ref{eq:esttl} implies the following existence of minimizer of the functional \;$J$.
\begin{cor}\label{eq:corexistencemin}
 Let \;$(M, g)$\; be a closed $n$-dimensional Riemannian manifold with \;$n\geq 4$\;even such that \;$\ker P_g^n\simeq \R$\; \;$P^n_g\geq 0$\; and \;$\kappa_g^n=(n-1)!\omega_n$. Assuming that \;$K$\; is a smooth positive function on \;$M$\; such that \;$(ND)_+$\; holds, then \;$K$\; is the \;$Q$-curvature of a Riemannian metric \;$g_u=e^{2u}g$\; with \;$u$\; a minimizer of \;$J$\; on \;$H^2(M)$.
\end{cor}

\noindent
\vspace{6pt}
 
\noindent
Theorem \ref{eq:esttl} implies also the following compactness theorem for conformal metrics with prescribed \;$Q$-curvature.
\begin{cor}\label{eq:compact}
Let \;$(M, g)$\; be a closed $n$-dimensional Riemannian manifold with \;$n\geq 4$\;even such that \;$\ker P_g^n\simeq \R$\; and \;$\kappa_g^n=(n-1)!m\omega_n$\; with $m\in \N^*$. Assuming that \;$K$\; is a smooth positive function on \;$M$\; such that \;$(ND)_0$\; holds, then for every \;$k\in \N$, there exists a large positive constant \;$C_k$\; such that for every \;$u$\; solution of \eqref{eq:qequation}, 
$$
||u||_{C^k(M)}\leq C_k.
$$
\end{cor}
\vspace{6pt}

\noindent
Corollary \ref{eq:compact} implies that the Leray-Schauder degree of equation \eqref{eq:qequation} is well-defined. Indeed, we have
\begin{thm}\label{eq:cordegree}
Let \;$(M, g)$\; be a closed $n$-dimensional Riemannian manifold with \;$n\geq 4$\;even such that \;$\ker P_g^n\simeq \R$\; and \;$\kappa_g^n=(n-1)!m\omega_n$\; with $m\in \N^*$. Assuming that $K$ is a smooth positive function on $M$ such that $(ND)_0$ holds, then the Leray-Schauder degree \;$d_m$\; of equation \eqref{eq:qequation} is well-defined. Furthermore, if \;$(ND)$\; holds, then  there existe $L_0>0$ such that
\begin{equation}
\begin{split}
 d_m=
\begin{cases}
(-1)^{\bar m}\left(1-\sum_{A\in \mathcal{F}_{\infty}} (-1)^{i_{\infty}(A)}\right)=\chi(J^L, J^{-L})-\sum_{A\in \mathcal{F}_{\infty}} (-1)^{\bar m+i_{\infty}(A)}\;\;&\text{if}\;\;m=1,\\
(-1)^{\bar m}\left(\frac{1}{(m-1)!}\Pi_{i=1}^{m-1}(i-\chi(M))- \frac{1}{m!}\sum_{A\in \mathcal{F}_{\infty}} (-1)^{i_{\infty}(A)}\right)\\=\chi(J^L, J^{-L})- \frac{1}{m!}\sum_{A\in \mathcal{F}_{\infty}} (-1)^{\bar m+i_{\infty}(A)}\;\;&\text{if}\;\;\;m\geq 2.
\end{cases}
\end{split}
\end{equation}
for all \;$L\geq L_0$.
\end{thm}

\section{Notations and preliminaries}\label{eq:notpre}
In this brief section, we fix our notations and give some preliminaries. First of all, we recall that \;$(M, \;g)$\,  and \;$K$\; are respectively the given underlying closed \;$n$-dimensional Riemannian manifold with \,$n\geq 4$\; even, and the prescribed function with the following properties (until otherwise said)
\begin{equation}\label{eq:trivialker}
\ker P_g^n\simeq \R\;\;\text{ and} \;\;\kappa_g^{n}=(n-1)!m\omega_n\;\; \text {for some }\;\;m\in \N^* \;\;\text{and}\;\;K\;\;\text{is a smooth positive function on} \;M.
\end{equation}
Furthermore, $u_l$ and $t_l$ are respectively sequence of smooth functions and real numbers, satisfying
\begin{equation}\label{eq:blowupeq1}
P_g^nu_l+t_lQ_g^n=t_l\kappa_g^nKe^{nu_l}\;\;\text{in}\;\;M,\;\;\text{and}\;\;t_l\longrightarrow 1\;\;\text{as}\;\;l\longrightarrow +\infty.
\end{equation}
\vspace{4pt}

\noindent
We are going to discuss the asymptotics near the singularity of the Green's function $G$ of the operator \;$P_{g}^n(\cdot)+\frac{1}{m}Q_g^n$\; with mass \;$(n-1)!\omega_n$\; satisfying the normalization $\int_MG(\cdot, y)Q^n_g(y)dV_g(y)=0$ \; and make some related definitions. 
\vspace{4pt}

\noindent
In the following, for a Riemmanian metric $\bar g$ on $M$, we will use the notation\;$B^{\bar g}_{p}(r)$\; to denote the geodesic ball with respect to $\bar g$ of radius \;$r$\;and center \;$p$. We also denote by \;$d_{\bar g}(x,y)$\; the geodesic distance with respect to $\bar g$ between two points \;$x$\;and \;$y$\; of \;$M$, $exp_x^{\bar g}$ the exponential map with respect to $\bar g$ at $x$. $inj_{\bar g}(M)$\;stands for the injectivity radius of \;$(M, \bar g)$, $dV_{\bar g}$\;denotes the Riemannian measure associated to the metric\;$\bar g$.  Furthermore, we recall that $\n_{\bar g}$, \;$\D_{\bar g}$, \;$R_{\bar g}$\; and \;$Ric_{\bar g}$\; will denote respectively the gradient, the Laplace-Beltrami operator, the scalar curvature and Ricci curvature with respect to \;$\bar g$. For simplicity, we will use the notation $B_p(r)$ to denote $B^g_{p}(r)$, namely $B_p(r)=B^g_p(r)$. $M^2$\;stands for the cartesian product \;$M\times M$, while \;$Diag(M)$\; is the diagonal of \;$M^2$.
\vspace{4pt}

\noindent
For \;$1\leq p\leq \infty$\; and \;$k\in \N$, \;$\theta\in  ]0, 1[$, \;$L^p(M)$, $W^{k, p}(M)$, $C^k(M)$, and $C^{k, \theta} (M)$ stand respectively for the standard Lebesgue space, Sobolev space, $k$-continuously differentiable space and $k$-continuously differential space of H\"older exponent $\beta$, all with respect $g$ (if the definition needs a metric structure, and for precise definitions and properties, see \cite{aubin} or \cite{gt}).  Given a function \;$u\in L^1(M)$,\;$\bar u$ and $\ov {u}_{Q^n}$\; denote respectively its average on \;$M$ with respect to $g$ and the sign measure $Q_g^ndV_g$, that is \;$$\bar u=\frac{\int_{M} u(x)dV_{g}(x)}{Vol_g(M)},$$ with \;$Vol_{g}(M)=\int_{M}dV_{g}$ and 
\begin{equation}\label{eq:wmass}
\ov{u}_{Q^n}=\frac{1}{(n-1)!\omega_n m}\int_{M} Q_g^n(x)u(x)dV_{g}(x).
\end{equation}
\vspace{4pt}

\noindent


\vspace{6pt}

 \noindent
 For  \;$\epsilon>0$\; and small,  $\l\in \R_+$, $\l\geq \frac{1}{\epsilon}$, and \;$a\in M$, $O_{\l, \epsilon}(1)$\; stands for quantities bounded uniformly in \;$\l$, and $\epsilon$, and \;$O_{a, \epsilon}(1)$ stands for quantities bounded uniformly in $a$ and $\epsilon$. For  $l\in \N^*$, $O_{l}(1)$ stands for quantities bounded uniformly in \;$l$\; and \;$o_l(1)$ stands for quantities which tends to $0$ as $l\rightarrow +\infty$.  For  $\epsilon$ positive and small, \;$a\in M$\; and \;$\l\in \R_+$ large, $\l\geq \frac{1}{\epsilon}$,\;$O_{a, \l, \epsilon}(1)$\; stands for quantities bounded uniformly in \;$a$, \;$\l$, and $\epsilon$. For $\epsilon$ positive and small, $p\in \N^{*}$, $\bar \l:=(\l_1, \cdots, \l_p)\in (\R_+)^p$, $\l_i\geq \frac{1}{\epsilon}$  for $i=1, \cdots, p$, and $A:=(a_1, \cdots, a_p)\in M^p$ (where $ (\R_+)^p$ and \;$M^p$\; denotes respectively the cartesian product of $p$ copies of $\R_+$ and $M$), $O_{A, \bar \l, \epsilon}(1)$ stands for quantities bounded uniformly in $A$, $\bar \l$, and $\epsilon$. Similarly for $\epsilon $ positive and small,  $p\in \N^{*}$, $\bar \l:=(\l_1, \cdots, \l_p)\in (\R_+)^p$, $\l_i\geq \frac{1}{\epsilon}$ for $i=1, \cdots, p$, $\bar \alpha:=(\alpha_1, \cdots, \alpha_p)\in \R^p$, $\alpha_i$ close to $1$ for $i=1, \cdots, p$, and $A:=(a_1, \cdots, a_p)\in M^p$ (where $ \R^p$  denotes the cartesian product of $p$ copies of $\R$, $O_{\bar\alpha, A, \bar \l, \epsilon}(1)$ will mean quantities bounded from above and below independent of $\bar \alpha$, $A$, $\bar \l$, and $\epsilon$. For $x\in \R$, we will use the notation $O(x)$ to mean $|x|O(1)$ where $O(1)$ will be specified in all the contexts where it is used. Large positive constants are usually denoted by $C$ and the value of\;$C$\;is allowed to vary from formula to formula and also within the same line. Similarly small positive constants are also denoted by $c$ and their value may varies from formula to formula and also within the same line.
\vspace{6pt}

\noindent
We call \;$\bar m$\; the number of negative eigenvalues (counted with multiplicity) of \;$P_g^n$. We point out that \;$\bar m$\; can be zero, but it is always finite. If \;$\bar m\geq 1$, then we will denote by \;$E_{-}\subset W^{\frac{n}{2}, 2}(M)$ the direct sum of the eigenspaces corresponding to the negative eigenvalues of \;$P_g^n$. The dimension of $E_{-}$ is of course $\bar m$. On the other hand, we have the existence of an \;$L^2$-orthonormal basis of eigenfunctions \;$v_1,\cdots, v_{\bar m}$\; of \;$E_{-}$ satisfying
\begin{equation}\label{eq:defeigen1}
P_g^n v_i=\mu_r v_r\;\;\;\forall\;\;r=1\cdots \bar m,
\end{equation}
\begin{equation}\label{eq:defeigen2}
\mu_1\leq \mu_2\leq \cdots\leq \mu_{\bar m}<0<\mu_{\bar m+1}\leq\cdots,
\end{equation}
where \;$\mu_r$'s are the eigenvalues of \;$P_g^n$\; counted with multiplicity. We define also the pseudo-differential operator \;$P^{n, +}_g$\; as follows
\begin{equation}\label{eq:operatorrev}
P_g^{n, +}u=P_g^nu-2\sum_{r=1}^{\bar m}\m_r(\int_M u v_rdV_g) v_r.
\end{equation}
Basically \;$P_g^{n, +}$\; is obtained from \;$P_g^n$ by reversing the sign of the negative eigenvalues and we extend the latter definition to \;$\bar m=0$ for uniformity in the analysis and recall that in that case \;$P_g^{n, +}=P_g^n$.  Now, 
 for \;$t>0$\; we set
\begin{equation}\label{eq:jt1}
J_t(u):=\left\langle P_g^{n}u, u\right\rangle+2t\int_MQ_g^nudV_g-t\frac{2(n-1)!\omega_nm}{n}\log\int_MKe^{nu}dV_g, \,\;\;u\in W^{\frac{n}{2}, 2}(M),
\end{equation}
and hence \;$J=J_1$. 
\vspace{6pt}

\noindent
We will use the notation \;$\left\langle\cdot, \cdot\right\rangle$\; to denote the \;$L^2$ scalar product and \;$\left\langle\cdot, \cdot\right\rangle_{W^{\frac{n}{2}, 2}}$\; for the \;$W^{\frac{n}{2}, 2}$-scalar product. On the other hand, it is easy to see that
\begin{equation}\label{eq:productpr}
\left\langle u, v\right\rangle_{P^n}:=\left\langle P^{n, +}_gu, v\right\rangle, \;\,\;u, v\in\{w\in W^{\frac{n}{2}, 2}(M):\;\;\;\ov u_{Q^n}=0\}
\end{equation}
defines a inner product on $\{u\in W^{\frac{n}{2}, 2}(M):\;\;\;u_{Q^n}=0\}$ which induces a norm equivalent to the standard  norm \;$||\cdot||:=\sqrt{\left\langle\cdot, \cdot\right\rangle_{W^{\frac{n}{2}, 2}}}$\; of \;$W^{\frac{n}{2}, 2}(M)$\;  (on $\{u\in W^{\frac{n}{2}, 2}(M):\;\;\;u_{Q^n}=0\}$) and denoted by
\begin{equation}\label{eq:normpr}
||u||_{P^n}:=\sqrt{\left\langle u, u\right\rangle_{P^n}} \;\;\;u\in\{w\in W^{\frac{n}{2}, 2}(M):\;\;\;\ov u_{Q^n}=0\}.
\end{equation}
$\bar B^{\bar m}_r$ will stand for the closed ball of center \;$0$ and radius \;$r$\; in \;$\R^{\bar m} $. $\mathbb{S}^{\bar m-1}$\; will denote the boundary of \;$\bar B^{\bar m}_1$. Given a set \;$X$,  we define \;$\widetilde{X\times \bar B^{\bar m}_1}$\; to be the cartesian product \;$X\times \bar B^{\bar m}_1$\; where the tilde means that \;$X\times \partial B^{\bar m}_1$\; is identified with \;$\partial B_1^{\bar m}$. For \;$m\geq 2$, we denote by \;$B_{m-1}(M)$\; the set of formal barycenters of \;$M$\; of order \;$m-1$, namely
 \begin{equation}\label{eq:barytop}
B_{m-1}(M):=\{\sum_{i=1}^{m-1}\alpha_i\d_{a_i}, a_i\in M, \alpha_i\geq 0, i=1,\cdots, m-1,\;\,\sum_{i=1}^{m-1}\alpha_i=1\},
\end{equation}
 Finally, we set
\begin{equation}\label{eq:defbarynega}
A_{m-1, \bar m}:=\widetilde{B_{m-1}(M)\times \bar B^{\bar m}_1}.
\end{equation}
\vspace{4pt}

\noindent
In the sequel also, \;$J^{c}$\; with $c\in \R$ will stand for \;$J^{c} :=\{u\in W^{\frac{n}{2}, 2}(M):\;\;J(u)\leq c\}$. For $X$ a topological space,  $\chi(X)$  denotes the Euler characteristic of $X$\; with \;$\Z_2$\; coefficients and for \;$(X, Y)$\; a topological pair, $\chi(X, Y)$\; denotes the Euler characteristic of $X$\; with \;$\Z_2$\; coefficients.

\vspace{4pt}

\noindent
As above, in the general case, namely \;$\bar m\geq 0$, for $\epsilon$\; small and positive, $\bar \beta:=(\beta_1, \cdots, \beta_{\bar m})\in \R^{\bar m}$\;  with \;$\beta_i$\; close to \;$0$, $i=1, \cdots, \bar m$) (where \;$\R^{\bar m}$\; is the empty set when \;$\bar m=0$), $\bar \l:=(\l_1, \cdots, \l_p)\in (\R_+)^p$, $\l_i\geq \frac{1}{\epsilon}$\;  for $i=1, \cdots, p$, $\bar \alpha:=(\alpha_1, \cdots, \alpha_p)\in \R^p$, $\alpha_i$ close to $1$ for $i=1, \cdots, p$, and $A:=(a_1, \cdots, a_p)\in M^p$, $p\in \N^{*}$, $w\in W^{\frac{n}{2}, 2}$\; with \;$||w||$\; small, $O_{\bar\alpha, A, \bar \l, \bar \beta, \epsilon}(1)$\; will stand quantities bounded independent of \;$\bar \alpha$, $A$, $\bar \l$, $\bar \beta$, and \;$\epsilon$, and $O_{\bar\alpha, A, \bar \l, \bar \beta, w, \epsilon}(1)$\; will stand quantities bounded independent of \;$\bar \alpha$, $A$, $\bar \l$, $\bar \beta$,  $w$\; and \;$\epsilon$.
\vspace{6pt}

\noindent
As in  \cite{no1} and \cite{nd4}, given a point \;$b\in \R^n$\; and \;$\lambda$\; a positive real number, we define \;$\delta_{b, \lambda}$\; to be the {\em standard bubble}, namely
\begin{equation}\label{eq:standarbubble}
\delta_{b, \lambda}(y):=\log\left(\frac{2\lambda}{1+\lambda^2|y-b|^2}\right),\;\;\;\;\;\;y\in \R^n.
\end{equation}
The functions \;$\delta_{b, \lambda}$\; verify the following equation
\begin{equation}\label{eq:bubbleequation}
(-\D_{\R^n})^{\frac{n}{2}}\delta_{b,\lambda}=(n-1)!e^{n\delta_{b,\lambda}}\;\;\;\text{in}\;\;\;\R^n.
\end{equation}
Geometrically, equation \eqref{eq:bubbleequation} means that the metric \;$g=e^{2\delta_{b, \lambda}} dx^2$\; (after pull-back by stereographic projection) has constant \;$Q$-curvature equal to $(n-1)!$, where $dx^2$ is the standard metric on \;$\R^n$.\\
Using the existence of conformal normal coordinates (see \cite{cao} or \cite{gun}), we have that, for \;$a \in M$\; there exists a function \;$u_a\in C^{\infty}(M)$ such that
\begin{equation}\label{eq:detga}
g_a = e^{2u_a} g\;\; \text{verifies}\;\;det g_a(x)=1\;\;\text{for}\;\;\; x\in B^{g_a}_a( \varrho_a).
\end{equation}
with \;$0<\varrho_0<\varrho_a<\frac{inj_{g_a}(M)}{10}$ for some small positive \;$\varrho_0$\; satisfying \;$\varrho_0<\frac{inj_g(M)}{10}$.
\vspace{6pt}

\noindent
Now, for \;$0<\varrho<\min\{\frac{inj_g(M)}{4}, \frac{\varrho_0}{4}\}$, we define a smooth cut-off function \;$\chi_{\varrho} : \bar\R_+ \rightarrow \bar\R_+$\; satisfying the following properties:
\begin{equation}\label{eq:cutoff}
\begin{cases}
\chi_{\varrho}(t)  = t \;\;&\text{ for } \;\;t \in [0,\varrho],\\
\chi_{\varrho}(t) = 2 \varrho \;\;&\text{ for } \;\; t \geq 2 \varrho, \\
 \chi_{\varrho}(t) \in [\varrho, 2 \varrho] \;\;\;&\text{ for } \;\; t\in [\varrho, 2 \varrho].
\end{cases}
\end{equation}
Using the cut-off function $\chi_{\varrho}$, we define for $a\in M$ and $\lambda\in \R_+$  the function $\hat{\delta}_{a, \lambda}$ as follows
\begin{equation}\label{eq:hatdelta}
\hat{\delta}_{a, \lambda}(x):=\log\left(\frac{2\lambda}{1+\lambda^2\chi_{\varrho}^2(d_{g_a}(x, a))}\right).
\end{equation}
For every \;$a\in M$\; and \;$\lambda\in \R_+$, we define \;$\varphi_{a, \lambda}$\; to be the unique the solution of the following projection problem
\begin{equation}\label{eq:projbubble}
\begin{cases}
P_g^n \varphi_{a,\l} \, + \, \frac{1}{m} Q_g^n\, = (n-1)!\omega_n\,\frac{ e^{n (\hat{\d}_{a,\l} \,  +  \, u_a)}}{\int_M e^{n (\hat{\d}_{a,\l} \,  +  \, u_a)}dV_g}\;\; \mbox{ in } \;\;M, \\
\int_M Q_g^n(x)\varphi_{a,\l}(x) \,dV_g(x) = 0.
\end{cases}
\end{equation}
Now, we recall that \;$G$\; is the unique solution of the following PDE
\begin{equation}\label{eq:defG4}
\begin{cases}
P_g^n G(a, \cdot)+\frac{1}{m}Q^n_g(\cdot)=(n-1)!\omega_n\delta_a(\cdot),\\
\int_M Q_g^n(x) G(a, x)dV_g(x)=0.
\end{cases}
\end{equation}
Using  \eqref{eq:defG4}, it is easy to see that the following integral representation formula holds
\begin{equation}\label{eq:G4integral}
u(x)-\ov{u}_{Q^n}=\frac{1}{(n-1)!\omega_n}\int_M G(x, y)P_g^nu(y), \;\;\;u\in C^n(M), \;x\in M,
\end{equation}
where $u_{Q^n}$ is defined as in section \ref{eq:notpre}. It is a well know fact that \;$G\;$ has a logarithmic singularity. In fact \;$G$\; decomposes as follows
\begin{equation}\label{eq:decompG4}
 G(a,x)=\log \left(\frac{1}{\chi_{\varrho}^2(d_{{g}_a}(a, x))}\right)+H(a, x).
\end{equation}
where \;$H$\; is the regular par of \;$G$. Furthermore, it is also a well-know fact that
\begin{equation}\label{eq:regH4}
G\in C^{\infty}(M^2\setminus Diag(M)),\;\;\\,\;\text{and} \;H\in C^{3, \beta}(M^2)\;\;\;\forall \beta\in (0, 1).
\end{equation}
Now, using \eqref{eq:limitfs} and \eqref{eq:partiallimit} combined with the symmetry of \;$H$, it is easy to see that
\begin{equation}\label{eq:relationderivative}
\frac{\partial \mathcal{F}_K(a_1, \cdots, a_m)}{\partial a_i}=\frac{2}{n}\frac{\n_g\mathcal{F}^{A}_i(a_i)}{\mathcal{F}^{A}_i(a_i)}, \;\;\;i=1, \cdots, m.
\end{equation}
Next, we set
\begin{equation}\label{eq:deflA}
l_K(A):=\sum_{i=1}^m\left(\frac{\D_{g_{a_i}} \mathcal{F}^{A}_i(a_i)}{(\mathcal{F}^{A}_i(a_i))^{\frac{n-2}{n}}}-\frac{n}{2(n-1)}R_g(a_i)(\mathcal{F}^{A}_i(a_i))^{\frac{2}{n}}\right),
\end{equation}
and have
\begin{equation}\label{eq:auxiindexa1}
l_K(A)=\frac{2n}{n-2}\mathcal{L}_K(A), \;\;\forall A\in Crit(\mathcal{F}_K).
\end{equation}

\section{An expansion of  $\nabla J_{t_l}$ at infinity}
In this section, we present a useful expansion of \;$\n J_{t_l}$\;  at infinity. In order to do that, we first fix $\L$ to  a large positive constant. Next, like in \cite{nd5} and \cite{no1} (see also \cite{nd4}), for $\epsilon$ and $\eta$ small positive real numbers, we first denote by \;$V(m, \epsilon, \eta)$ the $(m,\epsilon, \eta)$-{\em neighborhood of potential critical points at infinity}, namely
\begin{equation}\label{eq:ninfinity}
\begin{split}
V(m, \epsilon, \eta):=\{u\in W^{\frac{n}{2}, 2}(M):a_1, \cdots, a_m\in M, \;\;\l_1,\cdots, \l_m>0, \;\;||u-\ov{u}_{Q^n}-\sum_{i=1}^m\varphi_{a_i,\lambda_i}||+\\ ||\n^{W^{\frac{n}{2},2}} J(u)||=O\left( \sum_{i=1}^m\frac{1}{\l_i}\right)\;\;\;\lambda_i\geq \frac{1}{\epsilon},\;\;\frac{2}{\Lambda}\leq \frac{\l_i}{\l_j}\leq \frac{\Lambda}{2},\;\;\text{and}\;\;d_g(a_i, a_j)\geq 4\ov C\eta\;\;\text{for}\;i\neq j\},
\end{split}
\end{equation}
where \;$\ov C$\; is a large positive constant,  $\n^{W^{\frac{n}{2}, 2}}J$ is the gradient of \;$J$\; with respect to the \;$W^{\frac{n}{2}, 2}$-topology, $O(1):=O_{A, \bar \l, u, \epsilon}(1)$ meaning bounded uniformly in \;$\bar\l:=(\l_1, \cdots, \l_n)$, $A:=(a_1, \cdots, m)$, $u$, $\epsilon$. Next,  as in  \cite{no1} and \cite{nd4}, following the ideas of Bahri-Coron\cite{bc1}, we have that for \;$\eta$\; a small positive real number with \;$0<2\eta<\varrho$, there exists \;$\epsilon_0=\epsilon_0(\eta)>0$ such that
\begin{equation}\label{eq:mini}
\forall\;0<\epsilon\leq \epsilon_0,\;\;\forall u\in V(m, \epsilon, \eta), \text{the minimization problem }\min_{B_{\epsilon, \eta}}||u-\ov{u}_{Q^n}-\sum_{i=1}^m\alpha_i\varphi_{a_i, \l_i}-\sum_{r=1}^{\bar m}\beta_r(v_r-\ov{(v_r)}_{Q^n})||_{P^n}
\end{equation}
has a unique solution, up to permutations, where \;$B_{\epsilon, \eta}$\; is defined as follows
\begin{equation}
\begin{split}
{B_{\epsilon, \eta}:=\{(\bar\alpha, A, \bar \l, \bar \beta)\in \R^m\times M^m\times (0, +\infty)^m\times \R^{\bar m}}:|\alpha_i-1|\leq \epsilon, \l_i\geq \frac{1}{\epsilon}, i=1, \cdots, m, \\d_g(a_i, a_j)\geq 4\ov C\eta, i\neq j, |\beta_r|\leq R, r=1, \cdots, \bar m\}.
\end{split}
\end{equation}
 Moreover, using the solution of \eqref{eq:mini}, we have that every $u\in V(m, \epsilon, \eta)$ can be written as
\begin{equation}\label{eq:para}
u-\ov{u}_{Q^n}=\sum_{i=1}^m\alpha_i\varphi_{a_i, \l_i}+\sum_{r=1}^{\bar m}\beta_r(v_r-\ov{(v_r)}_{Q^n}+w,
\end{equation}
where $w$ verifies the following orthogonality conditions
\begin{equation}\label{eq:ortho}
\begin{split}
\left\langle Q_g^n, w\right\rangle=\left\langle\varphi_{a_i, \l_i}, w\right\rangle_{P^n}=\left\langle\frac{\partial\varphi_{a_i, \l_i}}{\partial \l_i}, w\right\rangle_{P^n}=\left\langle\frac{\partial\varphi_{a_i, \l_i}}{\partial a_i}, w\right\rangle_{P^n}=\left\langle v_r, w\right\rangle=0, i=1, \cdots, m,\\r=1, \cdots, \bar m
\end{split}
\end{equation}
and the estimate
\begin{equation}\label{eq:estwmin}
||w||=O\left(\sum_{i=1}^m\frac{1}{\l_i}\right),
\end{equation}
where here $O\left(1\right):=O_{\bar \alpha, A, \bar \l, \bar\beta , w, \epsilon}\left(1\right)$,  and for the meaning of \;$O_{\bar \alpha, A, \bar \l, \bar\beta , w, \epsilon}\left(1\right)$, see Section \ref{eq:notpre}. Furthermore, the concentration points $a_i$,  the masses $\alpha_i$, the concentrating parameters $\l_i$ and the negativity parameter $\beta_r$ in \eqref{eq:para} verify also
\begin{equation}\label{eq:afpara}
\begin{split}
d_g(a_i, a_j)\geq 4\ov C\eta,\;i\neq j=1, \cdots, m, \frac{1}{\L}\leq\frac{\l_i}{\l_j}\leq\L\;\;i, j=1, \cdots, m, \;\;\l_i\geq\frac{1}{\epsilon},\;\;\text{and}\\\;\;\;\sum_{r=1}^{\bar m}|\beta_r|+\sum_{i=1}^m|\alpha_i-1|\sqrt{\log \l_i}=O\left(\sum_{i=1}^m\frac{1}{\l_i}\right)
\end{split}
\end{equation}
with still $O\left(1\right)$ as in \eqref{eq:estwmin}. Using the neighborhood of potential critical points at infinity, we have the following Lemma. For its proof see Lemma 3.1 in \cite{nd5} and Proposition 3. 3 in \cite{no1}.

\begin{lem}\label{eq:escape}
Let $\epsilon$ and $\eta$ be small positive real numbers with $0<2\eta<\varrho$ where $\varrho$ is as in \eqref{eq:cutoff}. Assuming that $u_l$ is a sequence of blowing up critical point  of $J_{t_l}$ with $\ov{(u_l)}_{Q^n}=0, l\in \N$ and $t_l\rightarrow 1$ as $l\rightarrow +\infty$ , then  there exists $l_{\epsilon, \eta}$ a large positive integer such that for every $l\geq l_{\epsilon, \eta}$, we have $u_l\in V(m, \epsilon, \eta)$.
\end{lem}
\vspace{6pt}

\noindent
Now, we present some gradient estimates  for  \;$J_{t_l}$. We start with the following one.
\begin{lem}\label{eq:gradientlambdaest}
Assuming that $\eta$ is a small positive real number with \;$0<2\eta<\varrho$ where $\varrho$ is as in \eqref{eq:cutoff}, and $\epsilon\leq \epsilon_0$ where $\epsilon_0$ is as in \eqref{eq:mini}, then  for $a_i\in M$ concentration points,  $\alpha_i$ masses , $\l_i$ concentration parameters ($i=1,\cdots,m$) and $\beta_r$ negativity parameters ($r=1, \cdots, \bar m$) satisfying \eqref{eq:afpara}, we have that for $l$ large enough and for every $j=1,\cdots, m$, there holds

\begin{equation}\label{gradlambda}
\begin{split}
&\left\langle\n J_{t_l}\left(\sum_{i=1}^m\alpha_i\varphi_{a_i, \l_i}+\sum_{r=1}^{\bar m}\beta_r(v_r-\ov{(v_r)}_{Q^n})\right), \l_j\frac{\partial \varphi_{a_j, \l_j}}{\partial \l_j}\right\rangle= 2(n-1)!\omega_n\alpha_j\tau_j\\&-\frac{c_n^2(n-1)!\omega_n}{n\l_j^2}\left(\frac{\D_{g_{a_j}} \mathcal{F}^{A}_j(a_j)}{\mathcal{F}^{A}_j(a_j)}-\frac{n}{2(n-1)}R_g(a_j)\right)-\frac{2(n-1)!\omega_n}{(n-2)\l_j^2}\tau_j \D_{g_{a_j}}H(a_j, a_j)\\&-\frac{2(n-1)!\omega_n}{(n-2)\l_j^2}\sum_{i=1, i\neq j}^m\tau_i \D_{g_{a_j}}G(a_j, a_i)+ \frac{c_n^2(n-1)\omega_n}{n\l_j^2}\tau_j\left(\frac{\D_{g_{a_j}} \mathcal{F}^{A}_j(a_j)}{\mathcal{F}^{A}_j(a_j)}-\frac{n}{2(n-1)}R_g(a_j)\right)\\&+O\left(\sum_{i=1}^m|\alpha_i-1|^2+\sum_{r=1}^{\bar m}|\beta_r|^3+\sum_{i=1}^m\frac{1}{\l_i^3}\right),
\end{split}
\end{equation}
where \;$A:=(a_1, \cdots, a_m)$, $c_n^2$\; is a positive real number depending only on \;$n$, and for \;$i=1, \cdots, m$, $$\tau_i:=1-t_l\frac{m\gamma_i}{D}, \;\;\;\;\;D:=\int_M K(x)e^{n(\sum_{i=1}^m\alpha_i\varphi_{a_i, \l_i}(x)+\sum_{r=1}^{\bar m}\beta_r v_r(x))}dV_g(x),$$ with 
$$\gamma_i:=c_i^n\l_i^{2n\alpha_i-n}\mathcal{F}^{A}_i(a_i)\mathcal{G}_i(a_i),$$ where 
$$
c_i^n:=\int_{\R^n}\frac{1}{(1+|y|^2)^{n\alpha_i}}dy
$$
\begin{equation*}
 \begin{split}
\mathcal{G}_i(a_i):=e^{n((\alpha_i-1)H(a_i, a_i)+\sum_{j=1, j\neq i}^m(\alpha_j-1)G(a_j, a_i))}e^{\frac{n}{2(n-2)}\sum_{j=1, j\neq i}^m\frac{\alpha_j}{\l_j^2}\D_{g_{a_j}}G(a_j, a_i)}e^{\frac{n}{2(n-2)}\frac{\alpha_i}{\l_i^2}\D_{g_{a_i}}H(a_i, a_i)}\\\times e^{n\sum_{r=1}^{\bar m}\beta_rv_r(a_i)},
\end{split}
\end{equation*}
$c_n^2$\; is a positive real number depending only on \;$n$\; and for the meaning of \;$O_{\bar\alpha, A, \bar \l, \bar \beta, \epsilon}\left(1\right)$, see Section \ref{eq:notpre}.
  Furthermore, we have
  \begin{equation}\label{gradlambdasum}
\begin{split}
&\left\langle\n J_{t_l}\left(\sum_{i=1}^m\alpha_i\varphi_{a_i, \l_i}+\sum_{r=1}^{\bar m}\beta_r(v_r-\ov{(v_r)}_{Q^n})\right), \sum_{i=1}^m\frac{\l_i}{\alpha_i}\frac{\partial \varphi_{a_i, \l_i}}{\partial \l_i}\right\rangle=\\& -\sum_{i=1}^m\frac{c_n^2(n-1)!\gamma_n}{n\l_i^2}\left(\frac{\D_{g_{a_i}} \mathcal{F}^{A}_i(a_i)}{\mathcal{F}^{A}_i(a_i)}-\frac{n}{2(n-1)}R_g(a_i)\right)+\bar c_n(1-t_l)m\\&+O\left(\sum_{i=1}^m|\alpha_i-1|^2+\sum_{r=1}^{\bar m}|\beta_r|^3+\sum^{m}_{i=1}\tau^3_i+\sum_{i=1}^m\frac{1}{\l_i^3}\right),
\end{split}
\end{equation}
where $A$, $O\left(1\right)$, $c^2_n$, and $\tau_i$ ($i=1, \cdots, m$) are as above and \;$\bar c_n$\; is positive constant depending only on \;$n$.
\end{lem}
\begin{pf}
The proof follows the same strategy as in Lemma 5.1  and Corollary 5.2 in \cite{no1}.
\end{pf}
\begin{rem}\label{eq.deptauil}
We would like to remark that the $\tau_i$'s depends on $t_l$, but for the seek of simplicity in notations, we have decided to omit this dependency. 
\end{rem}

\vspace{6pt}
 
\noindent
Next, we present a gradient estimate for $\n J_{t_l}$ in the direction of the $\alpha_i$'s. Indeed, we have:
\begin{lem}\label{eq:gradientalpha} 
Assuming that $\eta$ is a small positive real number with $0<2\eta<\varrho$ where $\varrho$ is as in \eqref{eq:cutoff}, and $0<\epsilon\leq \epsilon_0$ where $\epsilon_0$ is as in \eqref{eq:mini}, then  for $a_i\in M$ concentration points,  $\alpha_i$ masses, $\l_i$ concentration parameters ($i=1,\cdots,m$), and $\beta_r$ negativity parameters ($r=1, \cdots, \bar m$) satisfying \eqref{eq:afpara}, we have that for $l$ large enough and for every $j=1, \cdots, m$, there holds

\begin{equation}\label{gradalpha}
\begin{split}
&\left\langle\n J_{t_l}\left(\sum_{i=1}^m\alpha_i\varphi_{a_i, \l_i}+\sum_{r=1}^{\bar m}\beta_r(v_r-\ov{(v_r)}_{Q^n})\right), \varphi_{a_j, \l_j}\right\rangle=\\&\left(2\log\l_j+H(a_j, a_j)-C_2^{n}\right)\frac{1}{\alpha_j}\left\langle\n J_{t_l}\left(\sum_{i=1}^m\alpha_i\varphi_{a_i, \l_i}+\sum_{r=1}^{\bar m}\beta_r(v_r-\ov{(v_r)}_{Q^n})\right), \l_j\frac{\partial \varphi_{a_j, \l_j}}{\partial \l_j}\right\rangle\\&+\sum_{i=1, i\neq j}^m G(a_j, a_i)\left\langle\n J_{t_l}\left(\sum_{i=1}^m\alpha_i\varphi_{a_i, \l_i}+\sum_{r=1}^{\bar m}\beta_r(v_r-\ov{(v_r)}_{Q^n})\right), \l_i\frac{\partial \varphi_{a_i, \l_i}}{\partial \l_i}\right\rangle\\&+4(n-1)!\omega_n(\alpha_j-1)\log\l_j+O\left(\log\l_j\left(\sum_{i=1}^m\frac{|\alpha_i-1|}{\log \l_i}+(\sum_{r=1}^{\bar m}|\beta_r|)(\sum_{i=1}^m\frac{1}{\log \l_i})+\sum_{i=1}^m\frac{1}{\l^2_i}\right)\right),
\end{split}
\end{equation}
where \;$O\left(1\right)$\; is as in Lemma \ref{eq:gradientlambdaest} and \;$C^n_2$\; is a constant depending only on \;$n$.
\end{lem}
\begin{pf}
It follows from the same arguments as in the proof of Lemma 5.3 in \cite{no1}.
\end{pf}
\vspace{6pt}

\noindent
Now, we derive a gradient estimate for $\n J_{t_l}$ with respect to the $a_i$'s. Precisely, we have:
\begin{lem}\label{eq:gradientaest}
Assuming that $\eta$ is a small positive real number with $0<2\eta<\varrho$ where $\varrho$ is as in \eqref{eq:cutoff}, and $0<\epsilon\leq \epsilon_0$ where $\epsilon_0$ is as in \eqref{eq:mini}, then  for $a_i\in M$ concentration points,  $\alpha_i$ masses, $\l_i$ concentration parameters ($i=1,\cdots,m$), ad $\beta_r$ negativity parameters ($r=1, \cdots, \bar m$) satisfying \eqref{eq:afpara}, we have that for $l$ large enough and for every $j=1, \cdots, m$, there holds
\begin{equation}\label{grada}
\begin{split}
\left\langle\n J_{t_l}\left(\sum_{i=1}^m\alpha_i\varphi_{a_i, \l_i}+\sum_{r=1}^{\bar m}\beta_r(v_r-\ov{(v_r)}_{Q^n})\right), \frac{1}{\l_j}\frac{\partial \varphi_{a_j, \l_j}}{\partial a_j}\right\rangle&=-\frac{4c^2_n(n-1)!\omega_n}{n\l_j}\frac{\n_g\mathcal{F}_j^{A}(a_j)}{\mathcal{F}_j^{A}(a_j)}\\&+O\left(\sum_{i=1}^m|\alpha_i-1|^2\right)\\&+O\left(\sum_{i=1}^m\frac{1}{\l_i^2}+\sum_{r=1}^{\bar m}|\beta_r|^2+\sum_{i=1}^m\tau_i^2\right),
\end{split}
\end{equation}
where $A:=(a_1, \cdots, a_m)$, $O(1)$ is as in Lemma \ref{eq:gradientlambdaest} and for  $i=1, \cdots, m$, \;$\tau_i$ is as in Lemma \ref{eq:gradientlambdaest}.
\end{lem}
\begin{pf}
It follows from the same arguments as in the proof of Lemma 5.4 in \cite{no1}.
\end{pf}
\vspace{6pt}

\noindent
Finally, we have the following estimate for \;$\n J_{t_l}$ in the direction of the $\beta_r$'s.
\begin{lem}\label{eq:gradientbeta}
Assuming that $\eta$ is a small positive real number with $0<2\eta<\varrho$ where $\varrho$ is as in \eqref{eq:cutoff}, and $0<\epsilon\leq \epsilon_0$ where $\epsilon_0$ is as in \eqref{eq:mini}, then for $a_i\in M$ concentration points,  $\alpha_i$ masses, $\l_i$ concentration parameters ($i=1,\cdots,m$), ad $\beta_r$ negativity parameters ($r=1, \cdots, \bar m$) satisfying \eqref{eq:afpara}, we have that for $l$ large enough and for every $s=1, \cdots, \bar m$, there holds
here holds
\begin{equation}\label{gradbeta}
\begin{split}
\left\langle\n J_{t_l}\left(\sum_{i=1}^m\alpha_i\varphi_{a_i, \l_i}+\sum_{r=1}^{\bar m}\beta_r (v_r-\ov{(v_r)}_{Q^n})\right), v_l-\ov{(v_l)}_Q\right\rangle=&2\mu_l\beta_l+O\left(\sum_{i=1}^m|\alpha_i-1|+\sum_{i=1}^m |\tau_i| +\sum_{i=1}^m\frac{1}{\l_i^2}\right),
\end{split}
\end{equation}
where $O(1)$ is as in Lemma \ref{eq:gradientlambdaest} and for  $i=1, \cdots, m$, \;$\tau_i$ is as in Lemma \ref{eq:gradientlambdaest}.
\end{lem}
\begin{pf}
 It follows from the same arguments as in the proof of Lemma 5.5 in \cite{no1}.
\end{pf}
\begin{rem}\label{remc1}
We would like to point out that the gradient estimates in \eqref{gradlambda}-\eqref{gradbeta} holds in \;$C^1$\; as function of the variables \;$(\bar \alpha, A, \bar \lambda, \bar \beta, \tau)$\; with $\bar \alpha=(\alpha_1, \cdots, \alpha_m), A=(a_1, \cdots, a_m), \bar \l=(\l_1, \cdots, \l_m), \bar \beta=(\beta_1, \cdots, \beta_{\bar m})$, $\tau=(\tau_1, \cdots, \tau_m)$  where for $B, C, D$\; some $C^1$-functions of the variables \;$(\bar \alpha, A, \bar \lambda, \bar \beta, \tau)$, $B=C+O(D)$ in $C^1$\; in the variables \;$(\bar \alpha, A, \bar \lambda, \bar \beta, \tau)$\; means \;$B=C+O(D)$\; and $\tilde \n B=\tilde \n C+O(D)$\; with \;$\tilde \n$\; denoting the gradient with respect to \;$(\bar \alpha, A, \bar \lambda, \bar \beta, \tau)$.
\end{rem}
\section{Refined location of $u_l$}\label{eq:reflocul}
In this section, we improve the location of $u_l$ given by Lemma \ref{eq:escape}. In order to do that, we divide this section into two subsections. In the first one, we derive a finite-dimensional parametrization of $u_l$. In the second one, we present the improvement we talked about above, by using the later parametrization of $u_l$.
\subsection{Finite-dimensional parametrization of $u_l$}
As already mentioned above, in this subsection, we give a finite-dimensional parametrization of $u_l$. For this end, we start with the following Lemma.
\begin{lem}\label{eq:expansionJ1}

Assuming that \;$\eta$\; is a small positive real number with \;$0<2\eta<\varrho$\; where \;$\varrho$\; is as in \eqref{eq:cutoff}, and \;$0<\epsilon\leq \epsilon_0$\; where \;$\epsilon_0$\; is as in \eqref{eq:mini} and \;$u=\ov{u}_Q+\sum_{i=1}^m\alpha_i\varphi_{a_i, \l_i}+\sum_{r=1}^{\bar m}\beta_r (v_r-\ov{(v_r)}_{Q^n})+w\in V(m, \epsilon, \eta)$ with \;$w$, the concentration points \;$a_i$,  the masses \;$\alpha_i$, the concentrating parameters \;$\l_i$  ($i=1, \cdots, m$), and the negativity parameters \;$\beta_r$ ($r=1, \cdots, \bar m$) verifying \eqref{eq:ortho}-\eqref{eq:afpara}, then we have
\begin{equation}\label{eq:exparoundbubble}
J_{t_l}(u)=J_{t_l}\left(\sum_{i=1}^m\alpha_i\varphi_{a_i, \l_i}+\sum_{r=1}^{\bar m}\beta_r (v_r-\ov{(v_r)}_{Q^n})\right)-f_l(w)+Q_l(w)+o(||w||^2),
\end{equation}
where
\begin{equation}\label{eq:linear}
f_l(w):=2(n-1)!\omega_nt_l\frac{\int_M Ke^{n\sum_{i=1}^m\alpha_i\varphi_{a_i, \l_i}+n\sum_{r=1}^{\bar m}\beta_r v_r}wdV_g}{\int_M Ke^{n\sum_{i=1}^m\alpha_i\varphi_{a_i, \l_i}+n\sum_{r=1}^{\bar m}\beta_r v_r}dV_g},
\end{equation}
and
\begin{equation}\label{eq:quadratic}
Q_l(w):=||w||_{P^n}^2-n!\omega_nt_l\frac{\int_M Ke^{n\sum_{i=1}^m\alpha_i\varphi_{a_i, \l_i}+n\sum_{r=1}^{\bar m}\beta_r v_r}w^2dV_g}{\int_M Ke^{n\sum_{i=1}^m\alpha_i\varphi_{a_i, \l_i}+n\sum_{r=1}^{\bar m}\beta_r v_r}dV_g}.
\end{equation}
Moreover, setting
\begin{equation}\label{eq:eali}
\begin{split}
E_{a_i, \l_i}:=\{w\in W^{\frac{n}{2}, 2}(M): \;\;\langle\varphi_{a_i, \l_i}, w\rangle_{P^n}=\langle\frac{\partial\varphi_{a_i, \l_i}}{\partial \l_i}, w\rangle_{P^n}=\langle\frac{\partial\varphi_{a_i, \l_i}}{\partial a_i}, w\rangle_{P^n}=0,\\\,\;\,\;\langle w, Q_g^n\rangle=\langle v_r, w\rangle=0, \;r=1, \cdots, \bar m,\,\;\text{and}\;\;||w||=O\left(\sum_{i=1}^m\frac{1}{\l_i}\right)\},
\end{split}
\end{equation}
and
\begin{equation}\label{eq:eal}
A:=(a_1, \cdots, a_m), \;\;\bar \l=(\l_1, \cdots, \l_m), \;\;E_{A, \bar \l}:=\cap_{i=1}^m E_{a_i, \l_i},
\end{equation}
we have that, the quadratic form \;$Q$\; is positive definite in \;$E_{A, \bar \l}$. Furthermore, the linear part \;$f$\; verifies that, for every \;$w\in E_{A, \bar \l}$, there holds
\begin{equation}\label{eq:estlinear}
f_l(w)=O\left[ ||w||\left(\sum_{i=1}^m\frac{|\n_g \mathcal{F}^{A}_i(a_i)|}{\l_i}+\sum_{i=1}^m|\alpha_i-1|\log \l_i+\sum_{r=1}^{\bar m}|\beta_r|+\sum_{i=1}^m\frac{\log \l_i}{\l_i^{2}}\right)\right].
\end{equation}
where here $o(1)=o_{ l, \bar\alpha, A, \bar \beta, \bar\l, w, \epsilon}(1)$ and $O\left(1\right):=O_{l, \bar\alpha, A, \bar \beta, \bar\l, w, \epsilon}\left(1\right)$ and for their meanings see section \ref{eq:notpre}.
\end{lem}
\begin{pf}
 The proof is the same as the one Proposition 6.1 in \cite{no1}.
\end{pf}

\vspace{6pt}

\noindent
Like in \cite{no1} and for the same reasons, we have that Lemma \ref{eq:expansionJ1} implies the following direct corollary.
\begin{cor}\label{eq:c1expansionj1}
Assuming that $\eta$ is a small positive real number with $0<2\eta<\varrho$ where $\varrho$ is as in \eqref{eq:cutoff}, $0<\epsilon\leq \epsilon_0$ where $\epsilon_0$ is as in \eqref{eq:mini}, and $u:=\sum_{i=1}^m\alpha_i\varphi_{a_i, \l_i}+\sum_{r=1}^{\bar m}\beta_r(v_r-\ov{(v_r)}_{Q^n})$ with the concentration  points $a_i$,  the masses $\alpha_i$, the concentrating parameters $\l_i$ ($i=1, \cdots, m$)  and the negativity parameters $\beta_r$ ($r=1, \cdots, \bar m$) satisfying \eqref{eq:afpara}, then for $l$ large enough, there exists a unique $\bar w_l(\bar \alpha, A, \bar \l, \bar \beta)\in E_{A, \bar \l}$ such that
\begin{equation}\label{eq:minj}
J_{t_l}(u+\bar w_l(\bar \alpha, A, \bar \l, \bar \beta))=\min_{w\in E_{A, \bar \l}, u+w\in V(m, \epsilon, \eta)} J_{t_l}(u+w),
\end{equation}
where $\bar \alpha:=(\alpha_1, \cdots, \alpha_m)$, $A:=(a_1, \cdots, a_m)$, $\bar \l:=(\l_1, \cdots, \l_m)$ and $\bar \beta:=(\beta_1, \cdots, \beta_m)$.\\
Furthermore, for $l$ large enough, the map $(\bar \alpha, A, \bar \l, \bar \beta)\longrightarrow \bar w_l(\bar \alpha, A, \bar \l, \bar \beta)\in C^1$ and satisfies the following estimate
\begin{equation}\label{eq:linminqua}
\frac{1}{C}||\bar w_l(\bar \alpha, A,  \bar \l, \bar \beta)||^2\leq |f_l(\bar w_l(\bar \alpha, A, \bar \l, \bar \beta))|\leq C|| \bar w_l(\bar \alpha, A, \bar \l, \bar \beta)||^2,
\end{equation}
for some large positive constant $C$ independent of $l$, $\bar \alpha$, $A$, $\bar \l$, and $\bar \beta$, hence
\begin{equation}\label{eq:estbarw}
||\bar w_l(\bar \alpha, A, \bar \l, \bar \beta)||=O\left(\sum_{i=1}^m\frac{|\n_g \mathcal{F}^{A}_i(a_i)|}{\l_i}+\sum_{i=1}^m|\alpha_i-1|\log\l_i+\sum_{r=1}^{\bar m}|\beta_r|+\sum_{i=1}^m\frac{\log \l_i}{\l_i^{2}}\right), 
\end{equation}
where $O\left(1\right):=O_{l, \bar\alpha, A, \bar \beta, \bar\l, \epsilon}\left(1\right)$ and for its meaning see section \ref{eq:notpre}. Moreover, assuming that
 $u_0:=\sum_{i=1}^m\alpha_i^0\varphi_{a_i^0, \l_i^0}+\sum_{r=1}^{\bar m}\beta_r^0(v_r-\ov{(v_r)}_{Q^n})$ with the concentration  points $a_i^0$,  the masses $\alpha_i^0$, the concentrating parameters $\l_i^0$ ($i=1, \cdots, m$) and the negativity parameters $\beta_r^0$ ($r=1, \cdots, \bar m$) satisfying \eqref{eq:afpara}, then for $l$ large enough, there exists an open neighborhood  $U^l$ of $(\bar\alpha^0, A^0, \bar \l^0, \bar \beta^0)$  (with $\bar \alpha^0:=(\alpha^0_1, \cdots, \alpha^0_m)$, $A^0:=(a_1^0, \cdots, a^0_m)$,  $\bar \l:=(\l_1^0, \cdots, \l_m^0)$ and  $\bar \beta^0:=(\beta_1^0, \cdots, \beta_{\bar m}^0)$) such that for every $(\bar \alpha, A, \bar \l, \bar \beta)\in U$ with $\bar \alpha:=(\alpha_1, \cdots, \alpha_m)$, $A:=(a_1, \cdots, a_m)$,  $\bar \l:=(\l_1, \cdots, \l_m)$, $\bar \beta:=(\beta_1, \cdots, \beta_{\bar m})$, and the $a_i$,  the $\alpha_i$, the $\l_i$ ($i=1, \cdots, m$)  and the $\beta_r$ ($r=1, \cdots, \bar m$) satisfying \eqref{eq:afpara}, and $w$ satisfying \eqref{eq:afpara}  with $\sum_{i=1}^m\alpha_i\varphi_{a_i, \l_i}+\sum_{r=1}^{\bar m}\beta_r(v_r-\ov{(v_r)}_{Q^n}+w\in V(m, \epsilon, \eta)$, we have the existence of a change of variable
\begin{equation}\label{eq:changev}
w\longrightarrow V_l
\end{equation}
from a neighborhood of $ \bar w_l(\bar \alpha, A, \bar \l, \bar \beta)$ to a neighborhood of $0$ such that
\begin{equation}\label{eq:expjv}
\begin{split}
&J_{t_l}(\sum_{i=1}^m\alpha_i\varphi_{a_i, \l_i}+\sum_{r=1}^{\bar m}\beta_r(v_r-\ov{(v_r)}_{Q6n})+w)=\\&J_{t_l}(\sum_{i=1}^m\alpha_i\varphi_{a_i, \l_i}+\sum_{r=1}^{\bar m}\beta_r(v_r-\ov{(v_r)}_{Q^n})+\bar w_l (\bar \alpha, A, \bar \l, \bar \beta))\\&+\frac{1}{2}\partial^2 J_{t_l}(\sum_{i=1}^m\alpha_i^0\varphi_{a_i^0, \l_i^0}+\sum_{r=1}^{\bar m}\beta_r^0(v_r-\ov{(v_r)}_{Q^n})+\bar w_l(\bar \alpha^0, A^0, \bar \l^0, \bar \beta^0))(V_l, V_l).
\end{split}
\end{equation}
\end{cor}
\vspace{6pt}

\noindent
Thus, with this new variable, in $V(m, \epsilon, \eta)$ we have a splitting of the variables $(\bar \alpha, A, \bar \l, \bar \beta)$ and $V_l$,  and $-V_l$  is a pseudogradient of $J_{t_l}$ in the direction of $V_l$. Using this fact, we have the following Proposition  which was the goal of this subsection.
\begin{pro}\label{eq:optimalul}
Assuming that $u_l$ is a sequence of blowing up solutions to \eqref{eq:blowupeq}, then for $l$ large enough there holds
\begin{equation}\label{eq:finpara}
u_l-\ov{(u_l)}_{Q^n}=\sum_{i=1}^m\alpha_i^l\varphi_{a_i^l, \l_i^l}+\sum_{r=1}^{\bar m}\beta_r^l(v_r-\ov{(v_r)}_{Q^n})+\bar w_l(\bar \alpha_l, A_l, \bar \l_l, \bar \beta_l),
\end{equation}
with $\bar\alpha_l:=(\alpha^l_1, \cdots, \l_m^l)$, $A_l:=(a_1^l, \cdots, a_m^l)$, $\bar\l_l:=(\l_1^l, \cdots, \l_m^l)$, and $\bar\beta_l=(\beta_1^l, \cdots, \beta_{\bar m}^l)$. Furthermore, for $l$ large enough, the concentration  points $a_i^l$, the masses $\alpha_i^l$, the concentrating parameters $\l_i^l$ ($i=1, \cdots, m$) and the negativity parameters $\beta_r^l$ ($r=1, \cdots, \bar m$) satisfy \eqref{eq:afpara}.
\end{pro}
\begin{pf}
It follows from the fact that $u_l$ is a solution to \eqref{eq:blowupeq} implies $\n J_{t_l}(u_l)=0$, the fact that $J_{t_l}$ is invariant by translation by constants combined with \eqref{eq:mini}, \eqref{eq:para}-\eqref{eq:afpara},  Lemma \ref{eq:escape}, Corollary \ref{eq:c1expansionj1} and the discussion right after it.
\end{pf}

\subsection{Refined estimates for the finite-dimensional parameters of \;$u_l$}
As already mentioned at the beginning of this section, in this subsection we derive refined estimates for the finite-dimensional parameters of $u_l$ in the formula \eqref{eq:finpara}. In order to do that, we start by constructing a pseudo-gradient for $J_l(\bar \alpha, A, \bar \l, \bar \beta)$, where $J_l(\bar \alpha, A, \bar \l, \bar \beta):=J_{t_l}(\sum_{i=1}^m\alpha_i\varphi_{a_i, \l_i}+\sum_{r=1}^{\bar m}\beta_r(v_r-\ov{(v_r)}_{Q^n})+\bar w_l (\bar \alpha, A, \bar \l, \bar \beta))$ in a suitable subset of $V(m, \epsilon, \eta)$. Indeed, setting
\begin{equation}\label{eq:ninfinityr}
\begin{split}
V_{deep}^{t_l}(m, \epsilon, \eta):= \{u\in V(m, \epsilon, \eta):\:\left(\sum_{i=1}^m\frac{|\n_g \mathcal{F}^{A}(a_i)|}{\l_i} +\sum_{i=1}^m|\alpha_i-1|+\sum_{i=1}^m|\tau_i|+\sum_{i=1}^m\frac{1}{\l_i^2}\right)\leq C_0\sum_{i=1}^m\frac{1}{\l_i^2}\},
 \end{split}
\end{equation}
with \;$C_0$\; a large positive constants, we have the following Proposition.
\begin{pro}\label{eq:conspseudograd}
Assuming that $\eta$ is a small positive real number with $0<2\eta<\varrho$ where $\varrho$ is as in \eqref{eq:cutoff}, and $0<\epsilon\leq \epsilon_0$ where $\epsilon_0$ is as in \eqref{eq:mini}, then for $l$ large enough, we have that there exists a pseudogradient $W:=W(l)$ of $J_l(\bar \alpha, A, \bar \l, \bar \beta)$  in \;$V(m, \epsilon, \eta)\setminus V_{deep}^{t_l}(m, \epsilon, \eta)$\; such that for every $u:=\sum_{i=1}^m\alpha_i\varphi_{a_i, \l_i}+\sum_{r=1}^{\bar m}\beta_r (v_r-\ov{(v_r)}_{Q^n})\in V(m, \epsilon, \eta)\setminus V_{deep}^{t_l}(m, \epsilon, \eta)$ with the concentration  points $a_i$, the masses $\alpha_i$, the concentrating parameters $\l_i$ ($i=1, \cdots, m$)  and the negativity parameters $\beta_r$ ($r=1, \cdots, \bar m$) satisfying \eqref{eq:afpara}, there holds
\begin{equation}\label{eq:pseudoexact}
<-\n J_l(u), W>\geq c\left(\sum_{i=1}^m\frac{1}{\l_i^2}+\sum_{i=1}^m\frac{|\n_g\mathcal{F}^{A}_i(a_i)|}{\l_i}+\sum_{i=1}^m|\alpha_i-1|+\sum_{i=1}^m|\tau_i|+\sum_{r=1}^{\bar m}|\beta_r|)\right),
\end{equation}
where $c$ is a small positive constant independent of $l$, $A:=(a_1, \cdots, a_m)$, $\bar\alpha=(\alpha_1, \cdots, \alpha_m)$, $\bar\l=(\l_1, \cdots, \l_m)$, $\bar \beta=(\beta_1, \cdots, \beta_{\bar m})$ and $\epsilon$. Furthermore, for every $u:=\sum_{i=1}^m\alpha_i\varphi_{a_i, \l_i}+\sum_{r=1}^{\bar m}\beta_r (v_r-\ov{(v_r)}_Q)+\bar w_l(\bar \alpha, A, \bar \l, \bar \beta)\in V(m, \epsilon, \eta)\setminus V_{deep}^{t_l}(m, \epsilon, \eta)$ with the concentration  points $a_i$, the masses $\alpha_i$, the concentrating parameters $\l_i$ ($i=1, \cdots, m$)  and the negativity parameters $\beta_r$ ($r=1, \cdots, \bar m$) satisfying \eqref{eq:afpara}, and $\bar w_l(\bar \alpha, A, \bar \l, \bar \beta)$ is as in \eqref{eq:minj},  there holds
\begin{equation}\label{eq:pseudoperturb}
<-\n J_l(u), W+\frac{\partial \bar w_l(W)}{\partial (\bar\alpha, A, \bar \l, \bar \beta)}>\geq c\left(\sum_{i=1}^m\frac{1}{\l_i^2}+\sum_{i=1}^m\frac{|\n_g\mathcal{F}^{A}_i(a_i)|}{\l_i}+\sum_{i=1}^m|\alpha_i-1|+\sum_{i=1}^m|\tau_i|+\sum_{r=1}^{\bar m}|\beta_r|\right),
\end{equation}
where \;$c$\; is still a small positive constant independent of $l$, $A:=(a_1, \cdots, a_m)$, $\bar\alpha=(\alpha_1, \cdots, \alpha_m)$, $\bar\l=(\l_1, \cdots, \l_m)$, $\bar \beta=(\beta_1, \cdots, \beta_{\bar m})$ and $\epsilon$.
\end{pro}
\begin{pf}
The argument is the same as the one of Proposition  in \cite{no1}. 
\end{pf}
 \vspace{6pt}
 
 \noindent
 Finally, we are going to achieve the goal of this section, by establishing a refined location of $u_l$, by exploiting its criticality for $J_{t_l}$, its finite-dimensional parametrization given by the previous subsection and Proposition \ref{eq:conspseudograd}. Precisely, we have:
\begin{lem}\label{eq:refloclem}
Let $\eta$ be a small positive real number with $0<2\eta<\varrho$ where $\varrho$ is as in \eqref{eq:cutoff} and $0<\epsilon\leq \epsilon_0$ where $\epsilon_0$ is as in \eqref{eq:mini}. Assuming that $u_l$ is a sequence of blowing-up solutions to \eqref{eq:blowupeq1}, then  for $l$ large enough, we have
$$
u_l\in V_{deep}^{t_l}(m, \epsilon, \eta).
$$
\end{lem}
\begin{pf}
It follows from the fact that \;$u_l$\; is a solution to \eqref{eq:blowupeq} implies $\n J_{t_l}(u_l)=0$ combined with Proposition \ref{eq:optimalul} and Proposition \ref{eq:conspseudograd}.
\end{pf}

\section{Proof of the results}\label{eq:pfresult}
In this section, we present the proof of Theorem \ref{eq:esttl}-Corollary \ref{eq:compact}. We start with the one of Theorem \ref{eq:esttl}. In order to do that, we are going to show the following result from which Theorem \ref{eq:esttl} follows directly, thanks to the formula \eqref{eq:auxiindexa1} and to the used scaling blowing up formula which is given by point a) of Lemma 2.3 in \cite{nd5}.
\begin{thm}\label{eq:esttla}
Let \;$(M, g)$\; be a closed four-dimensional Riemannian manifold such that \;$\ker P_g\simeq \R$, and  \;$\kappa_g^n=(n-1)!m\omega_n$\; with \;$m\in \N^*$. Assuming that $K$ is a smooth positive function on $M$, $\epsilon$ and $\eta$ be small positive real numbers with $0<2\eta<\varrho$ where $\varrho$ is as in \eqref{eq:cutoff} and $u_l$ is a sequence of blowing up solutions to \eqref{eq:blowupeq1}, then for $l$ large enough, we have that $u_l\in V_{deep}^{t_l}(m, \epsilon, \eta)$ 
$$
u_l-\ov{(u_l)}_{Q^n}=\sum_{i=1}^m\alpha_i^l\varphi_{a_i^l, \l_i^l}+\sum_{r=1}^{\bar m}\beta_r^l(v_r-\ov{(v_r)}_{Q^n})++\bar w_l(\bar \alpha_l, A_l, \bar \l_l, \bar \beta_l),
$$
and 
\begin{equation}\label{bubbleformulalambda}
t_l-1=\frac{\bar c_{n, m}^K(A^l)}{(\mathcal{F}^A(a_i))^{\frac{n-2}{n}}(\l_i^l)^2}\left[l_K(A^l)+O\left(\frac{1}{\l_i^l}\right)\right], \;\;i=1\cdots,m
\end{equation}
with $V_{deep}^{t_l}(m, \epsilon, \eta)$ defined by \eqref{eq:ninfinityr}, $\bar\alpha_l:=(\alpha^l_1, \cdots, \alpha_m^l)$, $A_l:=(a_1^l, \cdots, a_m^l)$, $\bar\l_l:=(\l_1^l, \cdots, \l_m^l)$,  $\bar\beta_l=(\beta_1^l, \cdots, \beta_{\bar m}^l)$,  the concentration  points $a_i^l$, the masses $\alpha_i^l$, the concentrating parameters $\l_i^l$ ($i=1, \cdots, m$), and the negativity parameters $\beta_r^l$ ($r=1, \cdots, \bar m$) satisfy \eqref{eq:afpara}, $A^l\longrightarrow A\in Crit(\mathcal{F}_K)$ as $l\rightarrow +\infty$, $\bar c_{n, m}^K(A^l)\longrightarrow \bar c_{n, m}^K(A)>0$, and $l_K(\cdot)$ is defined by \eqref{eq:deflA}.
\end{thm}
\begin{pfn}{ of Theorem \ref{eq:esttla}}\\
The proof is the same as the one of Theorem 5.1 in \cite{nd6}.
\end{pfn}
\vspace{6pt}

\noindent
\begin{pfn}{ of Corollary \ref{eq:corexistence}}\\
{\bf Case 1}: $(ND)_-$ \; holds\\
In this case, the result follows from our work \cite{nd1} in the {\em nonresonant} case by considering the functional \;$J_{1+\varepsilon}$ with $\varepsilon$ positive and small combined with Theorem \ref{eq:esttl}. \\\\
{\bf Case 2}: $(ND)_+$\; holds\\
In this case, the same argument as above works by considering \;$J_{1-\varepsilon}$\; with \;$\varepsilon$\; positive and small. We add that for \;$m=1$\; and \;$\bar m=0$, by the existence of minimizers in the subcritical case, we can take the solution to be be a minimizer of \;$J$.
On the other hand when \;$m\geq 2$, the result follows also from the characterization of the ``true'' critical points at infinity of  \;$J$, and the topology of very high and very negative sublevels of \;$J$\; established in our work  \cite{nd5}.
 \end{pfn}
\vspace{6pt}
 
\noindent

\begin{pfn}{ of Corollary \ref{eq:compact}}\\
Clearly Theorem \ref{eq:esttl} combined with standard elliptic regularity theory imply the result.
 \end{pfn}

\section{Proof of Theorem \ref{eq:cordegree}}
In this section, we prove Theorem \ref{eq:cordegree}. We start by showing a deep local characterization at infinity around critical points of \;$\mathcal{F}_K$\; for blowing-up solution of the type considered in Lemma \ref{eq:refloclem}. Indeed setting 
\begin{equation}\label{deepinfinity}
V_{deep}^{t_l}(m, \epsilon, \eta)(A^0):=\{u\in V_{deep}^{t_l}(m, \epsilon, \eta): \;\;\;d_g(a_i, a_i^0)\leq \tilde C_0\frac{1}{\l_i}, \;\;\;i;=1,\cdots, m\}
\end{equation}
for \;$A^0=(a_1^0, \cdots, a_m^0)\in Crit(\mathcal{F}_K)$\; with \;$\tilde C_0$\; a large positive constant, we have:
\begin{pro}\label{wdeep}
Let $\eta$ be a small positive real number with $0<2\eta<\varrho$ where $\varrho$ is as in \eqref{eq:cutoff} and $0<\epsilon\leq \epsilon_0$ where $\epsilon_0$ is as in \eqref{eq:mini}. Assuming that $u_l$ is a sequence of blowing-up solutions to \eqref{eq:blowupeq1}, then  for $l$ large enough, we have
$$
u_l\in V^{t_l}_{deep}(m, \epsilon, \eta)(A),
$$
for some \;$A\in Crit(\mathcal{F}_{K})$.
\end{pro}
\begin{pf}
Using Theorem \ref{eq:esttla}, we have \;$u_l\in V^{t_l}_{deep}(m, \epsilon, \eta)$\; and
$$
u_l-\ov{(u_l)}_{Q^n}=\sum_{i=1}^m\alpha_i^l\varphi_{a_i^l, \l_i^l}+\sum_{r=1}^{\bar m}\beta_r^l(v_r-\ov{(v_r)}_{Q^n})++\bar w_l(\bar \alpha_l, A_l, \bar \l_l, \bar \beta_l)
$$
with
$$
||\bar w_l(\bar \alpha^l, A, \bar \l^l, \bar \beta^l)||=O\left(\sum_{i=1}^m\frac{|\n_g \mathcal{F}^{A}_i(a_i^l)|}{\l_i^l}+\sum_{i=1}^m|\alpha_i^l-1|\log\l_i^l+\sum_{r=1}^{\bar m}|\beta_r^l|+\sum_{i=1}^m\frac{\log \l_i}{(\l_i^l)^{2}}\right).
$$
Thus using the definition of $V^{t_l}_{deep}(m, \epsilon, \eta)$, we have
$$
\frac{|\n_g \mathcal{F}^{A}_i(a_i^l)|}{\l_i^l}\leq C_0\frac{1}{(\l_i^l)^2}
$$
This implies \;$A^l=(a_i^l)\longrightarrow A\in Crit(\mathcal{F}_K)$, thank to \eqref{eq:relationderivative}. Moreover, the non-degeneracy of \;$\mathcal{F}_K$\; implies
$$
d_g(a_i^l, a_i)\leq \tilde C_0\frac{1}{\l_i^l}, \;\;i=1,\cdots, m,
$$
for some large \;$\tilde C_0>0$. Hence \;$u_l\in V^{t_l}_{deep}(m, \epsilon, \eta)(A)$. 
\end{pf}
\begin{rem}\label{vdeep1}
The bubbling rate formula in Theorem \ref{eq:esttl} implies that for \;$u_l\in V^{t_l}_{deep}(m, \epsilon, \eta)(A)$, we have  \;$u_l\in V^{1}_{deep}(m, \epsilon, \eta)(A)$, where $V^{1}_{deep}(m, \epsilon, \eta)(A)$ is defined as in \eqref{eq:ninfinityr} with \;$t_l$\; replaced by \;$1$\; and $A^0$\; replaced by \;$A$.
\end{rem}
\vspace{6pt}

\noindent
In the next proposition, we show that for any \;$A\in Crit(\mathcal{F}_K)$\; with \;$\mathcal{L}_K(A)<0$, there exists \;$u_t\in V^t_{deep}(m, \epsilon, \eta)(A)$\; for \;$t\simeq 1^{-}$\; when $\mathcal{F}_K$\; is a Morse function. The set \;$V^t_{deep}(m, \epsilon, \eta)(A)$\; is defined as in n \eqref{eq:ninfinityr} with \;$t_l$\; replaced by \;$t$\; and $A^0$\; replaced by \;$A$.
\begin{pro}\label{cons-bubble}
Let \;$(M, g)$\; be a closed \;$n$-dimensional Riemannian manifold with \;$n\geq 4$\;even such that \;$\ker P_g^n\simeq \R$\; and \;$\kappa_g^n=(n-1)!m\omega_n$. Assuming that \;$K$\; is a smooth positive function on \;$M$\; such that $\mathcal{F}_K$ is a Morse function and \;$A\in Crit(\mathcal{F}_K)$\; with \;$\mathcal{L}_K(A)<0$, then for \;$t\simeq 1^{-}$, there exist  \;$u_t$\; verifying \;
\begin{equation}\label{equationut}
P_g^nu_t+tQ^n_g=t\kappa^n_gKe^{nu_t}\;\;\text{in}\;\;M\end{equation} such that \;
\begin{equation}\label{neighborut}
u_t \in V^t_{deep}(m, \epsilon, \eta)(A).
\end{equation}
\end{pro}
\begin{pf}
By Theorem \ref{eq:esttla} and Proposition \ref{wdeep}, we must look for a solution \;$u^t\in V^t_{deep}(m ,\epsilon, \eta)(A)$\; for $t\simeq 1^{-}$\; verifying
$$
u_l-\ov{(u_t)}_{Q^n}=\sum_{i=1}^m\alpha_i^t\varphi_{a_i^t \l_i^t}+\sum_{r=1}^{\bar m}\beta_r^t(v_r-\ov{(v_r)}_{Q^n})+\bar w_t(\bar \alpha^t, A^t, \bar \l^t, \bar \beta^t),
$$
with  $\bar \alpha^t=(\alpha^t_1, \cdots, \alpha^t_m), \;A^t= (a^t_1, \cdots, a^t_m),\; \bar \l^t=( \l^t_1, \cdots, \l^t_m),\; \bar \beta^t=(\beta^t_1, \cdots, \beta^t_{\bar m})$,  and \;$\bar w_t(\bar \alpha^t, A^t, \bar \l^t, \bar \beta^t)$\; verifies $$
J_t(\sum_{i=1}^m\alpha_i^t\varphi_{a_i^t \l_i^t}+\sum_{r=1}^{\bar m}\beta_r^t(v_r-\ov{(v_r)}_{Q^n})+\bar w_t(\bar \alpha^t, A^t, \bar \l^t, \bar \beta^t))=\min_{w \in E_{A^t, \bar \l^t}, \hat z^t+w\in V(m, \epsilon, \eta)}J_t(\sum_{i=1}^m\alpha_i^t\varphi_{a_i^t \l_i^t}+\sum_{r=1}^{\bar m}\beta_r^t(v_r-\ov{(v_r)}_{Q^n})+w).
$$
with \;$$\hat z^t=\sum_{i=1}^m\alpha_i^t\varphi_{a_i^t \l_i^t}+\sum_{r=1}^{\bar m}\beta_r^t(v_r-\ov{(v_r)}_{Q^n})$$ and 
\begin{equation}\label{defwt}
||\bar w_l(\bar \alpha^t, A^t, \bar \l^t, \bar \beta^t)||=O\left(\sum_{i=1}^m\frac{|\n_g \mathcal{F}^{A^t}_i(a_i^t)|}{\l_i^t}+\sum_{i=1}^m|\alpha_i^t-1|\log\l_i^t+\sum_{r=1}^{\bar m}|\beta_r^t|+\sum_{i=1}^m\frac{\log \l_i^t}{(\l^t_i)^{2}}\right),
\end{equation}
where \;$E_{A^t, \bar \l^t}$\; is as in \eqref{eq:eal} with \;$(A, \bar \l)$\; replaced by \;$(A^t, \bar \l^t)$. 
Thus, we have 
$$
\left<\n J_t(\sum_{i=1}^m\alpha_i^t\varphi_{a_i^t \l_i^t}+\sum_{r=1}^{\bar m}\beta_r^t(v_r-\ov{(v_r)}_{Q^n})+\bar w_t(\bar \alpha^t, A^t, \bar \l^t, \bar \beta^t)), w\right>=0, \;\;\;\;\;\forall w\in E_{A^t, \bar \l^t}.
$$
Hence, setting \;
\begin{equation}\label{defzt}
z^t=\sum_{i=1}^m\alpha_i^t\varphi_{a_i^t \l_i^t}+\sum_{r=1}^{\bar m}\beta_r^t(v_r-\ov{(v_r)}_{Q^n})+\bar w_t(\bar \alpha^t, A^t, \bar \l^t, \bar \beta^t),
\end{equation}
we have 
\begin{equation}\label{fullsol}
\n J_t(z^t)=0
\end{equation}
is equivalent to
\begin{equation}\label{finitesol}
\begin{split}
\left<\n J_t(z^t),\l_j^t\frac{\partial \varphi_{a_j^t, \l_j^t} }{\partial \l_j^t}\right>=\left<\n J_t(z^t),\frac{1}{\l_j^t}\frac{\partial \varphi_{a_j^t, \l_j^t} }{\partial a_j^t}\right>=
\left<\n J_t(z^t), \varphi_{a_j^t, \l_j^t}\right>=\left<\n J_t(z^t), v_r-\ov{(v_r)}_{Q^n}\right>=0, \\ j=1, \cdots, m, \;r=1, \cdots \bar m.
\end{split}
\end{equation}
Using Lemma \ref{eq:expansionJ1} and Corollary \ref{eq:c1expansionj1}  with \;$t_l$\; replaced by \;$t$\; combined with \eqref{defwt} and recalling that we are looking for \;$z^t\in V^t_{deep}(m ,\epsilon, \eta)(A)$, we have
\begin{equation}\label{gradbetac1t}
\begin{split}
\left<\n J_t(z^t), v_r-\ov{(v_r)}_{Q^n}\right>=\left<\n J_t(\sum_{i=1}^m\alpha_i^t\varphi_{a_i^t \l_i^t}+\sum_{r=1}^{\bar m}\beta_r^t(v_r-\ov{(v_r)}_{Q^n})),  v_r-\ov{(v_r)}_{Q^n}\right>\\+O\left(\sum_{i=1}^m\left[\frac{|\n_g \mathcal{F}^{A^t}_i(a_i^t)|}{\l_i^t}\right]^2+\sum_{i=1}^m|\alpha_i^t-1|+\sum_{r=1}^{\bar m}|\beta_r^t|^2+\sum_{i=1}^n\frac{1}{(\l_i^t)^3}\right).
\end{split}
\end{equation}
As in Remark \ref{remc1}, \eqref{gradbetac1t} holds in $C^1$\; of the variables \;$(\bar \alpha^t, A^t, \bar \l^t, \bar \beta^t, \tau^t)$\; with $\tau^t=(\tau_1^t, \cdots, \tau_m^t)$\; with \;$\tau_i^t$\; as Lemma \ref{eq:gradientlambdaest} in with \;$t_l$ replaced by \;$t$. Thus, using Lemma \ref{eq:gradientbeta}, Remark \ref{remc1} and \eqref{gradbetac1t},  we have 
\begin{equation}\label{ztvr}
\left<\n J_t(z^t), v_r-\ov{(v_r)}_{Q^n}\right>=0
\end{equation}
is equivalent to
\begin{equation}\label{qequationbetat}
2\mu_r\beta_r^t+O\left(\sum_{i=1}^m\left[\frac{|\n_g \mathcal{F}^{A^t}_i(a_i^t)|}{\l_i^t}\right]^2+\sum_{i=1}^m|\alpha_i^t-1|+\sum_{i=1}^m|\tau_i^t|+\sum_{r=1}^{\bar m}|\beta_r^t|^2+\sum_{i=1}^m\frac{1}{(\l_i^t)^2}\right)=0,
\end{equation}
with \eqref{qequationbetat} holding in $C^1$ of the variables \;$(\bar \alpha^t, A^t, \bar \l^t, \bar \beta^t, \tau^t)$\; in the sense defined in Remark \ref{remc1}.
Thus by implicit function theorem we have a unique \;$\bar \beta^t_r=\bar \beta^t_r(A^t, \bar \alpha^t, \bar \l^t, \tau^t)$\; solving \eqref{qequationbetat} with \;$z^t\in V^t_{deep}(m ,\epsilon, \eta)(A)$\; for $t\simeq 1^{-}$. Similarly, using Lemma \ref{eq:gradientaest} and Corollary \ref{eq:c1expansionj1}  with \;$t_l$\; replaced by \;$t$\; combined with \eqref{defwt}, we have
\begin{equation}\label{gradat}
\begin{split}
\left<\n J_t(z^t), \frac{1}{\l_j^t}\frac{\partial \varphi_{a_j^t, \l_j^t} }{\partial a_j^t}\right>=\left<\n J_t(\sum_{i=1}^m\alpha_i^t\varphi_{a_i^t \l_i^t}+\sum_{r=1}^{\bar m}\beta_r^t(v_r-\ov{(v_r)}_{Q^n})), \frac{1}{\l_j^t}\frac{\partial \varphi_{a_j^t, \l_j^t} }{\partial a_j^t}\right>\\+O\left(\sum_{i=1}^m\left[\frac{|\n_g \mathcal{F}^{A^t}_i(a_i^t)|}{\l_i^t}\right]^2+\sum_{i=1}^m|\alpha_i^t-1|+\sum_{r=1}^{\bar m}|\beta_r^t|^2+\sum_{i=1}^n\frac{1}{(\l_i^t)^3}\right),
\end{split}
\end{equation}
with \;\eqref{gradat} holding in  \;$C^1$ of the variables \;$(\bar \alpha^t, A^t, \bar \l^t, \tau^t)$. Thus, using Lemma \ref{eq:gradientaest}, Remark \ref{remc1}, and \eqref{gradat}, we have 
\begin{equation}\label{zta}
\left<\n J_t(z^t), \frac{1}{\l_j^t}\frac{\partial \varphi_{a_j^t, \l_j^t} }{\partial a_j^t}\right>=0
\end{equation}\; is equivalent to
\begin{equation}\label{qequationat}
-\frac{4c^2_n(n-1)!\omega_n}{n\l_j}\frac{\n_g\mathcal{F}_j^{A^t}(a^t_j)}{\mathcal{F}_j^{A^t}(a^t_j)}+O\left(\sum_{i=1}^m\left[\frac{|\n_g \mathcal{F}^{A^t}_i(a_i^t)|}{\l_i^t}\right]^2+\sum_{i=1}^m|\alpha_i^t-1|+\sum_{i=1}^m\frac{1}{(\l_i^t)^2}+\sum_{r=1}^{\bar m}|\beta_r^t|^2+\sum_{i=1}^m(\tau_i^t)^2\right)=0,
\end{equation}
with \eqref{qequationat} holding in $C^1$ of the variables \;$(\bar \alpha^t, A^t, \bar \l^t, \tau^t)$. Thus as above, using the implicit function theorem we have a unique \;$A^t=A^t(\bar \alpha^t, \bar \l^t, \tau^t)$\; solving \eqref{qequationat} with \;$z^t\in V^t_{deep}(m ,\epsilon, \eta)(A)$\; for $t\simeq 1^{-}$. Again as above, using Lemma \ref{eq:expansionJ1} and Corollary \ref{eq:c1expansionj1}  with \;$t_l$\; replaced by \;$t$\; combined with \eqref{defwt}, we have
\begin{equation}\label{equationalphat}
\begin{split}
\left<\n J_t(z^t),  \varphi_{a_j^t, \l_j^t} \right>=\left<\n J_t(\sum_{i=1}^m\alpha_i^t\varphi_{a_i^t \l_i^t}+\sum_{r=1}^{\bar m}\beta_r^t(v_r-\ov{(v_r)}_{Q^n})), \varphi_{a_j^t, \l_j^t}\right>\\+O\left(\sum_{i=1}^m\left[\frac{|\n_g \mathcal{F}^{A^t}_i(a_i^t)|}{\l_i^t}\right]^2+\sum_{i=1}^m|\alpha_i^t-1|+\sum_{r=1}^{\bar m}|\beta_r^t|^2+\sum_{i=1}^n\frac{1}{(\l_i^t)^3}\right)
\end{split}
\end{equation}
with \eqref{equationalphat} holding in $C^1$ of the variables \;$(\bar \alpha^t,\bar \l^t, \tau^t)$. Thus, using Lemma \ref{eq:gradientlambdaest} and Lemma \ref{eq:gradientalpha}, Remark \ref{remc1}, and \eqref{equationalphat}, we have 
\begin{equation}
\begin{split}
\left<\n J_t(z^t),\varphi_{a_j^t, \l_j^t} \right>=0
\end{split}
\end{equation}
 is equivalent to
\begin{equation}\label{qequationalphat}
\begin{split}
4(n-1)!\omega_n(\alpha_j^t-1)\log\l_j^t+O\left(\sum_{i=1}^m\left[\frac{|\n_g \mathcal{F}^{A^t}_i(a_i^t)|}{\l_i^t}\right]^2+\sum_{i=1}^m|\alpha_i^t-1|+\sum_{i=1}^m\frac{1}{(\l_i^t)^2}+\sum_{r=1}^{\bar m}|\beta_r^t|+\sum_{i=1}^m(\tau_i^t)^2\right)=0.
\end{split}
\end{equation}
with \eqref{qequationalphat} holding in $C^1$ of the variables \;$(\bar \alpha^t,\bar \l^t, \tau^t)$. Thus, as above using the implicit function theorem we have a unique \;$\bar \alpha^t=\bar \alpha^t(\bar \l^t, \tau^t)$\; solving \eqref{qequationalphat} with \;$z^t\in V^t_{deep}(m ,\epsilon, \eta)(A)$\; for $t\simeq 1^{-}$. Again as above, using Lemma \ref{eq:expansionJ1} and Corollary \ref{eq:c1expansionj1}  with \;$t_l$\; replaced by \;$t$\; combined with \eqref{defwt},  we have
\begin{equation}\label{ztlambda}
\begin{split}
\left<\n J_t(z^t),  \l_j^t\frac{\partial \varphi_{a_j^t, \l_j^t} }{\partial \l_j^t}\right>=\left<\n J_t(\sum_{i=1}^m\alpha_i^t\varphi_{a_i^t \l_i^t}+\sum_{r=1}^{\bar m}\beta_r^t(v_r-\ov{(v_r)}_{Q^n})), \l_j^t\frac{\partial \varphi_{a_j^t, \l_j^t} }{\partial \l_j^t}\right>\\+O\left(\sum_{i=1}^m\left[\frac{|\n_g \mathcal{F}^{A^t}_i(a_i^t)|}{\l_i^t}\right]^2+\sum_{i=1}^m|\alpha_i^t-1|+\sum_{r=1}^{\bar m}|\beta_r^t|^2+\sum_{i=1}^n\frac{1}{(\l_i^t)^3}\right).
\end{split}
\end{equation}
with \eqref{ztlambda} holding in \;$C^1$\; of the variables \;$(\bar \l^t, \tau^t)$. Thus, using Lemma \ref{eq:gradientlambdaest}, Remark \ref{remc1} and \eqref{ztlambda}, we have 
\begin{equation}
\left<\n J_t(z^t),\l_j^t\frac{\partial \varphi_{a_j^t, \l_j^t} }{\partial \l_j^t}\right>=0
\end{equation}
is equivalent to
\begin{equation}\label{qequationlambdat}
2(n-1)!\omega_n\alpha_j^t\tau_j^t+O\left(\sum_{i=1}^m\left[\frac{|\n_g \mathcal{F}^{A^t}_i(a_i^t)|}{\l_i^t}\right]^2+\sum_{i=1}^m|\alpha_i^t-1|+\sum_{r=1}^{\bar m}|\beta_r^t|^2+\sum_{i=1}^n\frac{1}{(\l_i^t)^2}+\sum_{i=1}^m|\tau_i^t|^2\right)
\end{equation}
with \eqref{qequationlambdat} holding in $C^1$ of the variables \;$(\bar \l^t, \tau^t)$. Thus as above, using the implicit function theorem we have a unique \;$\tau^t=\tau^t(\bar \l^t)$\; solving \eqref{qequationlambdat} with \;$z^t\in V^t_{deep}(m ,\epsilon, \eta)(A)$\; for $t\simeq 1^{-}$. Finally for $t\simeq 1^{-1}$, we have $\forall i=1,\cdots, m$\; there exists a unique \;$\l^t_i$\; such that  
$$(t-1)(\mathcal{F}^{A}_i(a_i))^{\frac{n-2}{n}}=\frac{2n\bar c_{n m}^K(A)\mathcal{L}_K(A)}{n-2(\l_i^t)^2},$$
with \;$\bar c_{n , m}^K(A)$\; as in Theorem \ref{eq:esttla}. Hence 
$$
z^t=\sum_{i=1}^m\tilde \alpha_i^t\varphi_{\tilde a_i^t, \tilde  \l_i^t}+\sum_{r=1}^{\bar m}\tilde \beta_r^t(v_r-\ov{(v_r)}_{Q^n})+\bar w_t(\tilde \alpha^t, \tilde A^t, \tilde  \l^t, \tilde\beta^t)
$$
with $\tilde \l^t=\bar \l^t$, $\tilde \alpha^t=\bar \alpha^t(\bar \l^t, \tau^t(\bar \l^t))$, $\tilde A^t=A^t(\tilde \alpha^t, \bar \l^t, \tau^t(\bar \l^t))$, and $\tilde \beta^t=\bar \beta^t(\tilde \alpha^t, \bar\l^t,\tau^t(\bar \l^t) )$\; verifies
$$
\n J_t(z^t)=0.
$$
Thus \;$u^t=z^t-\log \int_Me^{4z^t}dV_g$\; satisfies \eqref{equationut} and \eqref{neighborut} thereby ending the proof of the proposition.
\end{pf}
\vspace{6pt}

\noindent
Our work  \cite{nd1} in the non-resonant case, Corollary \ref{eq:compact} and Theorem \ref{wdeep} imply the following proposition important in the calculation of \;$d_m$\; as carried below.
\begin{pro}\label{asymptotic}
Let \;$(M, g)$\; be a closed \;$n$-dimensional Riemannian manifold with \;$n\geq 4$\;even such that \;$\ker P_g^n\simeq \R$\; and \;$\kappa_g^n=(n-1)!m\omega_n$. Assuming that \;$K$\; is a smooth positive function on \;$M$\; such that \;$(ND)$\; holds. There exist \;$\epsilon_{m, n}>0$, $C_m>0$,  $C_t>C_{n, m}$ ($t\simeq 1$) with \;$C_{n, m}$\; and \;$\epsilon_{m, n}$ depending only on \;$m$\; and \;$n$, $C_t$\; continuous in \;$t$ and  \;$\lim _{t\rightarrow 1}C_t=+\infty$\; such that for every \;$u$\; solution of
$$
P^n_gu+tQ^n_g=tKe^{nu},
$$
with  \;$|(n-1)\omega_n(t-1)|<\epsilon_{m, n}$, we have\\
1)\\
$||u||< C_t\, \;\;\;\forall t\neq 1.$\\
2)\\ 
$||u||\leq C_{m, n}\, \;\;\text{for}\;\;t=1$\\
3) \\If $t\neq 1$, then we have\\
i) Either \;$||u||<C_{n ,m}$\\
ii) Or \;$||u||\geq C_{n ,m}$\; and  \;$u\in V^t_{deep}(m,\epsilon, \eta):=\cup_{A\in Crit(\mathcal{F}_K)}V^t_{deep}(m,\epsilon, \eta)(A)$.
\end{pro}
\begin{pf}
1) follows from the compactness result in \cite{nd1}( see also \cite{dr}, \cite{mal}), while 2) follows from Corollary \ref{eq:compact}, and 3) follows from Theorem \ref{wdeep}.
\end{pf}

\begin{pfn} {of Theorem \ref{eq:cordegree}}\\
Let\;$X:=\{u\in W^{\frac{n}{2}}(M):\;\;\int_MKe^{nu}dV_h=1\}$\; and \;$T: X\longrightarrow X$\; be defined by
\begin{equation}\label{opT}
T(u)=(P_g^n)^{-1}(\kappa_g^nKe^{nu}-Q^n_g), \;\;\;u\in X.
\end{equation}
Then \;$u$\; is a solution of \eqref{eq:qequation} is equivalent to \;$(I-T)u=0$. On the other hand, Corollary \ref{eq:compact} and Theorem \ref{asymptotic} imply the Leray-Schauder degree of \eqref{eq:qequation} \;$d_m=d_m(K)$\; is well-defined and verifies
$$
d_m=deg(I-T, B_{C_{n,m}}, 0).
$$
From the work of Malchiodi\cite{maldeg} in the non-resonant case, there exists \;$L_0>0$\; such that for all \;$L\geq L_0$,
$$
1-\chi(A_{m-1, \bar m})=\chi(J_t^{L}, J^{-L}_t), \;\;\;\forall t\in (0, 1).
$$
From our work \cite{nd5}, up to taking \;$L_0$\; larger we have 
$$
\chi(J^{L}, J^{-L})=\chi(J_t^{L}, J^{-L}_t), \;\;\;\forall t \in (0, 1).
$$
Similarly to \eqref{eq:qequation} and \eqref{opT}, for \;$t\in (0, 1)$\; we consider the equation 
\begin{equation}\label{qequationt}
P_g^nu+tQ^n_g=t\kappa^n_gKe^{nu},
\end{equation}
and the operator
$$
T_t(u)=(P_g^n)^{-1}(t\kappa_g^nKe^{nu}-tQ^n_g), \;\;\;u\in X.
$$
Then Proposition \ref{asymptotic} implies the Leray-Schauder degree of  \eqref{qequationt} \;$d_m^t$\; is well-defined and is given by
$$
d_m^t=deg (I-T_t, B_{C_t}, 0), \;\;\forall t\in (0, 1).
$$
Furthermore, the work of Malchodi\cite{maldeg} implies
$$
d_m^t=1-\chi(A_{m-1, \bar m}),\;\;\forall t\in (0, 1)
$$
Let use define
$$
d_m^{-}=\lim_{t\rightarrow 1^{-}}d_m^t.
$$
Then Theorem \ref{eq:esttl}, Remark \ref{vdeep1}, and Theorem \ref{asymptotic} imply
$$
d_m^{-}=deg(I-T, B_{C_{n, m}}, 0)+ deg(I-T, V^{1,-}_{deep}(m, \epsilon, \eta), 0),
$$
where
$$
V^{1,-}_{deep}(m, \epsilon, \eta)=\cup_{A\in \mathcal{F}_{\infty}}V^1_{deep}(m,\epsilon, \eta)(A)
$$
On the other hand, our Morse lemma at infinity in \cite{nd5}  and Poncare-Hopf theorem imply 
$$
deg(I-T,V^{1,-}_{deep}(m, \epsilon, \eta) , 0)=\frac{1}{m!}\sum_{A\in \mathcal{F}_{\infty}}(-1)^{i_{\infty}(A))+\bar m}
$$
Thus, we get
$$
d_m^{-}=d_m+\frac{1}{m!}\sum_{A\in \mathcal{F}_{\infty}}(-1)^{i_{\infty}(A)+\bar m}
$$
so, we obtain
$$
1-\chi(A_{m-1, \bar m})=\chi(J^{L}, J^{-L})=d_m+\frac{1}{m!}\sum_{A\in \mathcal{F}_{\infty}}(-1)^{i_{\infty}(A)}
$$
this implies
$$
d_m=\chi(J^{L}, J^{-L})-\frac{1}{m!}\sum_{A\in \mathcal{F}_{\infty}}(-1)^{i_{\infty}(A)+\bar m}=1-\chi(A_{m-1, \bar m})-\frac{1}{m!}\sum_{A\in \mathcal{F}_{\infty}}(-1)^{i_{\infty}(A)+\bar m}
$$
Hence, recalling
$$1-\chi(A_{m-1, \bar m})=(-1)^{\bar m}\;\;\text{for}\;\;\;m=1$$ and $$1-\chi(A_{m-1, \bar m})=(-1)^{\bar m}\frac{1}{(m-1)!}\Pi_{i=1}^{m-1}(i-\chi(M)),\;\;\text{for}\;\;\;m\geq2$$ 
we have the result follows.
\end{pfn}

\end{document}